\documentclass[3p]{elsarticle}

\usepackage[utf8]{inputenc}
\usepackage{amsmath,amsfonts,amssymb,amsthm,color}
\usepackage{graphicx}
\usepackage[labelsep=endash,labelfont=bf,justification=justified]{caption}
\usepackage[pdftex]{hyperref}
\usepackage{ulem} 
\usepackage{lineno}
\modulolinenumbers[5]
\allowdisplaybreaks

\begin{document}

\begin{frontmatter}

\title{Optimal control of PDEs using physics-informed neural networks}

\author[1,2]{Saviz Mowlavi}
\author[2]{Saleh Nabi\corref{cor1}}
\cortext[cor1]{Corresponding author}
\ead{nabi@merl.com}

\address[1]{Department of Mechanical Engineering, Massachusetts Institute of Technology, Cambridge, MA 02139, USA}
\address[2]{Mitsubishi Electric Research Laboratories, Cambridge, MA 02139, USA}

\begin{abstract}
Physics-informed neural networks (PINNs) have recently become a popular method for solving forward and inverse problems governed by partial differential equations (PDEs). By incorporating the residual of the PDE into the loss function of a neural network-based surrogate model for the unknown state, PINNs can seamlessly blend measurement data with physical constraints. Here, we extend this framework to PDE-constrained optimal control problems, for which the governing PDE is fully known and the goal is to find a control variable that minimizes a desired cost objective. We provide a set of guidelines for obtaining a good optimal control solution; first by \textcolor{black}{selecting an appropriate PINN architecture and training parameters based on a forward problem, second by choosing the best value for a critical scalar weight in the loss function using a simple but effective two-step line search strategy}. We then validate the performance of the PINN framework by comparing it to adjoint-based nonlinear optimal control, which performs gradient descent on the discretized control variable while satisfying the discretized PDE. This comparison is carried out on several distributed control examples based on the Laplace, Burgers, Kuramoto-Sivashinsky, and Navier-Stokes equations. Finally, we discuss the advantages and caveats of using the PINN and adjoint-based approaches for solving optimal control problems constrained by nonlinear PDEs.
\end{abstract}

\begin{keyword}
optimal control \sep physics-informed neural networks \sep adjoint-based optimisation \sep PDE-constrained optimization
\end{keyword}

\end{frontmatter}


\section{Introduction}

Owing to increases in computational power and wealth of available data as well as advances in algorithms, machine learning and deep learning specifically have revolutionized a number of fields over the last decade \cite{lecun2015deep} – image and speech recognition \cite{krizhevsky2012imagenet,hinton2012deep}, natural language processing \cite{sutskever2014sequence}, and drug discovery \cite{vamathevan2019applications} are just a few examples. In the physical sciences, however, data is more scarce while physical models are often available in the form of partial differential equations (PDEs) \cite{karniadakis2021physics}. 
Leveraging these governing equations, physics-informed neural networks (PINNs) were recently proposed in \cite{raissi2019physics} as a deep learning framework for solving forward and inverse problems without requiring much data (or not at all, in the case of forward problems). Building on a line of research that originated in the nineties \cite{dissanayake1994neural,van1995neural,lagaris1998artificial,hayati2007feedforward} and is rooted in the universal approximation theorem for neural networks \cite{hornik1989multilayer,leshno1993multilayer}, the basic idea behind PINNs is to approximate the solution to a given problem with a feed-forward neural network. This neural network is then trained by minimizing a composite loss function that not only penalizes the prediction error with respect to the available data but also enforces the governing equations and boundary conditions -- in effect, the governing equations act as an implicit prior that regularizes the training procedure in the data-limited regime. A great benefit of PINNs is their flexibility: they can either solve forward problems in the absence of any data when the governing equations are fully known, or leverage available data to solve inverse problems involving unknown model parameters or physical quantities (for reviews on PINNs, see \cite{lu2021deepxde} and \cite{cai2021physics}). PINNs have since found applications in numerous fields such as fluid mechanics \cite{raissi2019deep,raissi2020hidden,sun2020surrogate}, heat transfer \cite{cai2021physics}, solid mechanics \cite{rao2021physics,haghighat2021physics}, medicine \cite{sahli2020physics,van2020physics}, and chemistry \cite{ji2021stiff}, and they have been extended to account for noisy data \cite{yang2021b} as well as various types of governing equations such as stochastic PDEs \cite{zhang2020learning} and fractional PDEs \cite{pang2019fpinns}.

In this paper, we investigate the potential of PINNs to solve PDE-constrained optimal control problems, for which the governing PDEs are fully known and the goal is to find a control variable that minimizes a desired cost objective. Such problems arise in a variety of fields including fluid mechanics \cite{foures2014optimal}, transition to turbulence \cite{kerswell2018nonlinear}, heat transfer \cite{nabi2019nonlinear}, electromagnetism \cite{deng2018self}, topology optimization \cite{oktay2011parallelized}, and mesh refinement \cite{li2004adjoint}. The control variable to optimize might represent a distributed boundary actuation, an external body force or an initial condition of the system \cite{troltzsch2010optimal,borzi2011computational}. 
The highly nonlinear and multi-scale nature of such problems requires sophisticated numerical tools to determine control inputs that yield optimal performance according to application-specific criteria.
Optimal control problems are usually solved by combining gradient-descent algorithms with adjoint-based sensitivity analysis, which computes the gradient of the cost objective function with respect to the control variable using only two PDE simulations \cite{lions1971optimal}. This adjoint-based optimization framework is therefore very efficient when the control is a space- and/or time-dependent field, but its complexity has limited its adoption by the engineering community. By contrast, a major strength of PINNs is their ease of implementation, and we show here that the PINN framework can be readily extended to the optimal control setting by approximating the control field with its own neural network in addition to the neural network for the unknown state variable. These two networks are then simultaneously trained using a composite loss function that includes the cost objective functional in addition to the PDE residual and initial/boundary conditions. In this way, the training process finds a control and a state that satisfy the PDE constraint while minimizing the cost objective. A similar approach was recently proposed by \cite{lu2021physics} and \cite{demo2021extended} in the context of inverse design and parametric optimal control, respectively. In light of these recent works, the novelties of the present study are two-fold:
\begin{enumerate}
\item We propose a set of guidelines \textcolor{black}{consisting of two parts} for obtaining a good optimal control solution using the PINN framework. \textcolor{black}{First, we select an appropriate neural network architecture and training parameters by solving a forward problem based on the same PDE. Second, we solve the optimal control problem using a simple but effective line search strategy to find the best value for the scalar weight of the cost objective functional inside the loss function, which is a critical parameter to obtain a feasible optimal solution that simultaneously satisfies the PDE constraint while minimizing the cost objective. Specifically, we find the best scalar weight by} evaluating the cost objective \textcolor{black}{corresponding to the PINN optimal control} using a separate \textcolor{black}{PINN} forward computation that takes the PINN optimal control as input.
\item We apply these guidelines to solving a range of optimal control problems based on the Laplace, Burgers, Kuramoto-Sivashinsky, and Navier-Stokes equations using the PINN framework. In each case, we compare the quality of the optimal control found by PINNs with corresponding results obtained from adjoint-based optimization, which we implement on all examples. Finally, we leverage our experience from this comparative study to discuss the pros and cons of the PINN- and adjoint-based approaches for solving PDE-constrained optimal control problems. Such a careful and systematic comparison of the two methods enables researchers to better evaluate and position PINN-based optimal control within the larger context of PDE-constrained optimization.  
\end{enumerate}

This paper is structured as follows. The methodology and guidelines for solving optimal control problems with PINNs are presented in Section \ref{sec:Methodology}, along with a review of the classical adjoint-based optimization framework. The PINNs and adjoint-based approaches are then applied to a range of optimal control problems in Section \ref{sec:Results}, following which their pros and cons are discussed in Section \ref{sec:Discussion}. Conclusions close the paper in Section \ref{sec:Conclusions}.

\section{Methodology}
\label{sec:Methodology}

\subsection{Optimal control problem statement}

Consider a physical system defined over a domain $\Omega \subset \mathbb{R}^d$, and governed by a \textcolor{black}{PDE} of the form
\begin{subequations}
\begin{align}
\mathcal{F}[\mathbf{u}(\mathbf{x},t);\mathbf{c}_v(\mathbf{x},t)] &= 0, \quad \mathbf{x} \in \Omega, t \in [0,T], \label{eq:residual} \\
\mathcal{B}[\mathbf{u}(\mathbf{x},t);\mathbf{c}_b(\mathbf{x},t)] &= 0, \quad \mathbf{x} \in \partial \Omega, t \in [0,T], \label{eq:BC} \\
\mathcal{I}[\mathbf{u}(\mathbf{x},0);\mathbf{c}_0(\mathbf{x})] &= 0, \quad \mathbf{x} \in \Omega, \label{eq:IC}
\end{align} \label{eq:PDE}%
\end{subequations}
where $\mathbf{x}$ and $t$ denote respectively space and time, $\mathbf{u} \in \mathcal{U}$ is the state vector, $\mathbf{c}_v$, $\mathbf{c}_b$ and $\mathbf{c}_0$ are respectively volume, boundary and initial control vectors, $\mathcal{F}$ is the PDE residual which contains several differential operators, $\mathcal{B}$ are the boundary conditions, and $\mathcal{I}$ is the initial condition. The solution of the PDE \eqref{eq:PDE} depends on the control variable $\mathbf{c} = (\mathbf{c}_v,\mathbf{c}_b,\mathbf{c}_0) \in \mathcal{C}$. $\mathcal{U}$ and $\mathcal{C}$ are appropriate Hilbert spaces. For the discretized solution of PDEs, we assume such spaces to be, respectively, Euclidean spaces $\mathbb{R}^n$ and $\mathbb{R}^p$. In PDE-constrained optimal control, we seek the optimal $\mathbf{c}^*$ that minimizes a user-defined cost objective functional $\mathcal{J}(\mathbf{u}, \mathbf{c})$. This constrained optimization problem can be formulated as
\begin{equation}
\mathbf{c}^* = \arg \min_{\mathbf{c}} \mathcal{J}(\mathbf{u}, \mathbf{c}) \quad \text{subject to } \eqref{eq:PDE}.
\label{eq:OptimalControlProblem}
\end{equation}
A trivial approach to solve the above optimization problem is to calculate solutions of \eqref{eq:PDE} for a number of different control variables $\mathbf{c}$, and pick the one that leads to the lowest value of $\mathcal{J}$. However, when $\mathbf{c}$ is a continuous field defined over a finite volume and/or boundary, the resulting infinite dimensionality of the search space makes this approach impractical. We present in Section \ref{sec:PINNs} a methodology to solve the optimal control problem \eqref{eq:OptimalControlProblem} based on PINNs, followed in Section \ref{sec:AdjointBasedOptimalControl} by the classical framework of adjoint-based optimization which will provide our benchmark solutions in Section \ref{sec:Results}.

\textcolor{black}{Before proceeding further, we discuss the dependence of the cost objective functional on the control variable $\mathbf{c}$ in \eqref{eq:OptimalControlProblem}. In general, the goal of the problem is to find a control that results in a certain behavior for the state $\mathbf{u}$, which can be formulated, e.g., as a quadratic objective of the form $\mathcal{J}(\mathbf{u})$. However, the nonlinearity of the governing PDE can lead to the existence of multiple optimal control solutions $\mathbf{c}^*$. To guarantee the uniqueness of the optimal control solution, it is common practice to include a regularization term that depends on $\mathbf{c}$, resulting in a cost objective of the form $\mathcal{J}(\mathbf{u}, \mathbf{c})$. The regularization term can also be motivated by smoothness requirements for the control $\mathbf{c}$. For the three nonlinear PDEs considered in this paper, we included such a regularization term in the Kuramoto-Sivashinski example but we deliberately chose not to do so in the Burgers and Navier-Stokes examples for two reasons. First, one of our aims in this study is to compare the performance of the PINN and adjoint-based approaches at finding the best possible optimal control measured from a performance standpoint (i.e., how close the state $\mathbf{u}$ is to the desired behavior), which would be penalized by the presence of any regularization. Second, deep neural networks are empirically known to possess intrinsic regularization properties; a separate objective of this study is therefore to assess whether such inherent regularization would be reflected in the shape of the PINN optimal solution versus the adjoint-based solution, in the absence of explicit regularization. This being said, our problem formulation \eqref{eq:OptimalControlProblem} can be appplied to any choice of regularization.}

\subsection{Physics-informed neural networks for optimal control}
\label{sec:PINNs}

Let us first consider the forward problem defined by the PDE \eqref{eq:PDE} with prescribed control variable $\mathbf{c}$. In the framework of physics-informed neural networks, the solution $\mathbf{u}(\mathbf{x},t)$ of this forward problem is represented by a surrogate model $\mathbf{u}_{\mathrm{NN}}(\mathbf{x},t;\boldsymbol{\theta}_\mathbf{u})$ in the form of a fully-connected neural network that takes $(\mathbf{x},t)$ as input and returns an approximation for $\mathbf{u}$ at this $\mathbf{x}$ and $t$ as output. The vector $\boldsymbol{\theta}_\mathbf{u}$ denotes the set of trainable parameters of the network, which propagates the input data through its $\ell$ layers according to the sequence of operations
\begin{subequations}
\begin{align}
\mathbf{z}^0 &= (\mathbf{x},t), \label{eq:NNInput} \\
\mathbf{z}^k &= \sigma(\mathbf{W}^k \mathbf{z}^{k-1} + \mathbf{b}^k), \quad 1 \le k \le \ell-1, \\
\mathbf{z}^\ell &= \mathbf{W}^\ell \mathbf{z}^{\ell-1} + \mathbf{b}^\ell.
\end{align}
\end{subequations}
Each layer outputs a vector $\mathbf{z}^k \in \mathbb{R}^{q_k}$, where $q_k$ is the number of `neurons', and is defined by a weight matrix $\mathbf{W}^k \in \mathbb{R}^{q_k \times q_{k-1}}$, a bias vector $\mathbf{b}^k \in \mathbb{R}^{q_k}$, and a nonlinear activation function $\sigma(\cdot)$. Finally, the output of the last layer is used to represent the solution, that is, $\mathbf{u}_{\mathrm{NN}}(\mathbf{x},t;\boldsymbol{\theta}_\mathbf{u}) = \mathbf{z}^\ell$. \textcolor{black}{For prescribed} control variables $\mathbf{c} = (\mathbf{c}_v,\mathbf{c}_b,\mathbf{c}_0)$, the network parameters $\boldsymbol{\theta}_\mathbf{u} = \{\mathbf{W}^k,\mathbf{b}\}_{k=1}^\ell$ are trained by minimizing the loss function
\begin{align}
\mathcal{L}(\boldsymbol{\theta}_\mathbf{u},\mathbf{c}) &= \frac{1}{N_r} \sum_{i=1}^{N_r} |\mathcal{F}[\mathbf{u}_{\mathrm{NN}}(\mathbf{x}_i^r,t_i^r;\boldsymbol{\theta}_\mathbf{u});\mathbf{c}_v]|^2 + \frac{w_b}{N_b} \sum_{i=1}^{N_b} |\mathcal{B}[\mathbf{u}_{\mathrm{NN}}(\mathbf{x}_i^b,t_i^b;\boldsymbol{\theta}_\mathbf{u});\mathbf{c}_b]|^2 \nonumber \\
&\quad + \frac{w_0}{N_0} \sum_{i=1}^{N_0} |\mathcal{I}[\mathbf{u}_{\mathrm{NN}}(\mathbf{x}_i^0,0;\boldsymbol{\theta}_\mathbf{u});\mathbf{c}_0]|^2,
\label{eq:Loss}
\end{align}
where $\{\mathbf{x}_i^r, t_i^r\}_{i=1}^{N_r}$, $\{\mathbf{x}_i^b, t_i^b\}_{i=1}^{N_b}$, $\{\mathbf{x}_i^0\}_{i=1}^{N_0}$ each represent an arbitrary number of training points over which to enforce the PDE residual \eqref{eq:residual}, boundary conditions \eqref{eq:BC}, and initial condition \eqref{eq:IC}, respectively, and $w_b$, $w_0$ are scalar weights for the \textcolor{black}{boundary and initial loss components}. A critical underpinning of PINNs is the use of automatic differentiation (AD) to compute the loss \eqref{eq:Loss}. By using the chain rule to compose the derivatives of successive algebraic operations, AD calculates the exact derivatives of the network output $\mathbf{u}_{\mathrm{NN}}(\mathbf{x},t;\boldsymbol{\theta}_\mathbf{u})$ with respect to its inputs $\mathbf{x}$ and $t$. Thus, the various loss components in \eqref{eq:Loss} can be computed exactly without inheriting the truncation error incurred by standard numerical discretization schemes. Another advantage of computing derivatives with AD is that the residual points $\{\mathbf{x}_i, t_i\}_{i=1}^{N_r}$ can be chosen arbitrarily, conferring PINNs their convenient mesh-free nature. Starting from randomly initialized parameters $\boldsymbol{\theta}_\mathbf{u}$, we can now use gradient-based optimization to find an optimum set of values $\boldsymbol{\theta}_\mathbf{u}^*$ that minimizes \eqref{eq:Loss}. At each iteration $k$, the parameters are updated as
\begin{equation}
\boldsymbol{\theta}_\mathbf{u}^{k+1} = \boldsymbol{\theta}_\mathbf{u}^k - \alpha(k) \nabla_{\boldsymbol{\theta}_\mathbf{u}} \mathcal{L}(\boldsymbol{\theta}_\mathbf{u}^k;\mathbf{c}),
\label{eq:GradientUpdate}
\end{equation}
where $\alpha(k)$ is an adaptive learning rate set by the chosen optimizer. At the end of the training procedure, the trained neural network $\mathbf{u}_{\mathrm{NN}}(\mathbf{x},t;\boldsymbol{\theta}_\mathbf{u}^*)$ approximately solves the forward problem \eqref{eq:PDE}.

The PINN framework can be readily extended to the optimal control problem \eqref{eq:OptimalControlProblem}. Our approach is conceptually similar to the methodology proposed in \cite{lu2021physics} for solving PDE-constrained inverse design problems with PINNs, in which case $\mathbf{c}$ represents a set of design variables. For clarity of exposure, let us assume that we only have volume control, that is, $\mathbf{c} = \mathbf{c}_v$ (the treatment is similar for boundary or initial control). In addition to the neural network for $\mathbf{u}(\mathbf{x},t)$, we construct a second fully-connected neural network to approximate the control variable $\mathbf{c}(\mathbf{x},t)$. The neural network for $\mathbf{c}$ is denoted by $\mathbf{c}_{\mathrm{NN}}(\mathbf{x},t;\boldsymbol{\theta}_\mathbf{c})$, with $\boldsymbol{\theta}_\mathbf{c}$ the corresponding set of trainable parameters. Since we desire a solution for $\mathbf{u}$ and $\mathbf{c}$ that minimizes the cost \textcolor{black}{functional} $\mathcal{J}(\mathbf{u}, \mathbf{c})$ in addition to solving the PDE \eqref{eq:PDE}, we simply add \textcolor{black}{a term $\mathcal{L}_\mathcal{J}$ penalizing} the cost functional to the standard PINN loss \eqref{eq:Loss}, leading to the augmented loss function
\begin{align}
\mathcal{L}(\boldsymbol{\theta}_\mathbf{u}, \boldsymbol{\theta}_\mathbf{c}) &= \frac{w_r}{N_r} \sum_{i=1}^{N_r} |\mathcal{F}[\mathbf{u}_{\mathrm{NN}}(\mathbf{x}_i^r,t_i^r;\boldsymbol{\theta}_\mathbf{u});\mathbf{c}_{\mathrm{NN}}(\mathbf{x}_i^r,t_i^r;\boldsymbol{\theta}_\mathbf{c})]|^2 + \frac{w_b}{N_b} \sum_{i=1}^{N_b} |\mathcal{B}[\mathbf{u}_{\mathrm{NN}}(\mathbf{x}_i^b,t_i^b;\boldsymbol{\theta}_\mathbf{u})]|^2 \nonumber \\ &\quad + \frac{w_0}{N_0} \sum_{i=1}^{N_0} |\mathcal{I}[\mathbf{u}_{\mathrm{NN}}(\mathbf{x}_i^0,0;\boldsymbol{\theta}_\mathbf{u})]|^2 + w_\mathcal{J} \mathcal{L}_\mathcal{J}(\boldsymbol{\theta}_\mathbf{u}, \boldsymbol{\theta}_\mathbf{c}),
\label{eq:LossControl}
\end{align}
where we have introduced a new scalar weight $w_\mathcal{J}$ for the cost functional \textcolor{black}{term}. Due to the dependence of the PDE residual \eqref{eq:residual} on the control $\mathbf{c}$, the first loss term is now a function of both $\boldsymbol{\theta}_\mathbf{u}$ and $\boldsymbol{\theta}_\mathbf{c}$. The calculation of the term $\mathcal{L}_\mathcal{J}(\boldsymbol{\theta}_\mathbf{u}, \boldsymbol{\theta}_\mathbf{c})$ depends on the form on the objective functional and may involve residual points $\{\mathbf{x}_i^r, t_i^r\}_{i=1}^{N_r}$ or boundary points $\{\mathbf{x}_i^b, t_i^b\}_{i=1}^{N_b}$. Starting from randomly initialized parameters $(\boldsymbol{\theta}_\mathbf{u},\boldsymbol{\theta}_\mathbf{c})$, we can now use gradient-based optimization to find an optimum set of values $(\boldsymbol{\theta}_\mathbf{u}^*,\boldsymbol{\theta}_\mathbf{c}^*)$ that minimizes \eqref{eq:LossControl}. At each iteration $k$, the parameters from both networks are concurrently updated as
\begin{subequations}
\begin{align}
\boldsymbol{\theta}_\mathbf{u}^{k+1} &= \boldsymbol{\theta}_\mathbf{u}^k - \alpha(k) \nabla_{\boldsymbol{\theta}_\mathbf{u}} \mathcal{L}(\boldsymbol{\theta}_\mathbf{u}^k, \boldsymbol{\theta}_\mathbf{c}^k), \\
\boldsymbol{\theta}_\mathbf{c}^{k+1} &= \boldsymbol{\theta}_\mathbf{c}^k - \alpha(k) \nabla_{\boldsymbol{\theta}_\mathbf{c}} \mathcal{L}(\boldsymbol{\theta}_\mathbf{u}^k, \boldsymbol{\theta}_\mathbf{c}^k).
\end{align} \label{eq:GradientUpdateControl}%
\end{subequations}
At the end of the training procedure, the trained neural networks $\mathbf{u}_{\mathrm{NN}}(\mathbf{x},t;\boldsymbol{\theta}_\mathbf{u}^*)$ and $\mathbf{c}_{\mathrm{NN}}(\mathbf{x},t;\boldsymbol{\theta}_\mathbf{c}^*)$ approximately solve the optimal control problem \eqref{eq:OptimalControlProblem}.

In all results to follow, we use Glorot initialization of the parameters \cite{glorot2010understanding}, select a $\tanh$ activation function, and employ the Adam optimizer \cite{kingma2014adam}. We also normalize the input $(\mathbf{x},t)$ before passing it to the first layer of the neural network for $\mathbf{u}$ (and $\mathbf{c})$, so that \eqref{eq:NNInput} becomes
\begin{equation}
\mathbf{z}^0 = \left( \frac{\mathbf{x}-\mu_\mathbf{x}}{\sigma_\mathbf{x}}, \frac{t-\mu_t}{\sigma_t} \right),
\end{equation}
where $\mu_\mathbf{x}$, $\sigma_\mathbf{x}$, $\mu_t$, and $\sigma_t$ are the mean and standard deviation of the residual training points $\{\mathbf{x}_i, t_i\}_{i=1}^{N_r}$.

\subsection{Adjoint-based optimal control}
\label{sec:AdjointBasedOptimalControl}

The framework of adjoint-based optimal control is a direct extension of the method of Lagrange multipliers for constrained optimization to the case where the equality constraints are formulated as PDEs \cite{lions1971optimal}. Applying this method to the constrained problem \eqref{eq:OptimalControlProblem}, one first constructs the Lagrangian
\begin{equation}
\mathcal{L}(\mathbf{u},\mathbf{c},\boldsymbol{\lambda}) = \mathcal{J}(\mathbf{u},\mathbf{c}) - \langle \boldsymbol{\lambda}, \mathcal{F}[\mathbf{u};\mathbf{c}] \rangle,
\label{eq:Lagrangian}
\end{equation}
where $\mathbf{u}$ is required to satisfy the boundary and initial conditions \eqref{eq:BC} and \eqref{eq:IC}, $\boldsymbol{\lambda} = \boldsymbol{\lambda}(\mathbf{x},t)$ is the Lagrange multiplier or adjoint field, and the inner product $\langle \cdot, \cdot \rangle$ is defined as
\begin{equation}
\langle \mathbf{a}, \mathbf{b} \rangle = \int_0^T \int_\Omega \mathbf{a}(\mathbf{x},t)^\mathsf{T} \mathbf{b}(\mathbf{x},t) d\mathbf{x} dt.
\end{equation}
Then, the constrained problem \eqref{eq:OptimalControlProblem} is equivalent to the unconstrained problem
\begin{equation}
\mathbf{u}^*, \mathbf{c}^*, \boldsymbol{\lambda}^* = \arg \min_{\mathbf{u},\mathbf{c},\boldsymbol{\lambda}} \mathcal{L}(\mathbf{u}, \mathbf{c},\boldsymbol{\lambda}),
\label{eq:LagrangianProblem}
\end{equation}
whose solution is given by the stationary point(s) of the Lagrangian. This yields the relations
\begin{subequations}
\begin{align}
\left \langle \frac{\partial \mathcal{L}}{\partial \mathbf{u}}, \delta \mathbf{u} \right \rangle &= 0 \quad \forall \, \delta \mathbf{u}, \\
\left \langle \frac{\partial \mathcal{L}}{\partial \mathbf{c}}, \delta \mathbf{c} \right \rangle &= 0 \quad \forall \, \delta \mathbf{c}, \\
\left \langle \frac{\partial \mathcal{L}}{\partial \boldsymbol{\lambda}}, \delta \boldsymbol{\lambda} \right \rangle &= 0 \quad \forall \, \delta \boldsymbol{\lambda},
\end{align} \label{eq:StationarityConditions}%
\end{subequations}
where the admissible variation $\mathbf{u} + \delta \mathbf{u}$ has to satisfy the boundary and initial conditions \eqref{eq:BC} and \eqref{eq:IC}. The Fr\'echet derivative $\langle \partial \mathcal{L} / \partial \mathbf{u}, \cdot \rangle$ is defined so that
\begin{equation}
\left \langle \frac{\partial \mathcal{L}}{\partial \mathbf{u}}, \delta \mathbf{u} \right \rangle = \lim_{\epsilon \rightarrow 0} \frac{\mathcal{L}(\mathbf{u} + \epsilon \delta \mathbf{u},\mathbf{c},\boldsymbol{\lambda}) - \mathcal{L}(\mathbf{u},\mathbf{c},\boldsymbol{\lambda})}{\epsilon} \quad \forall \, \delta \mathbf{u},
\end{equation}
and similarly for $\langle \partial \mathcal{L} / \partial \mathbf{c}, \cdot \rangle$ and $\langle \partial \mathcal{L} / \partial \boldsymbol{\lambda}, \cdot \rangle$. Expanding the stationarity conditions \eqref{eq:StationarityConditions} leads to
\begin{subequations}
\begin{align}
\left \langle \frac{\partial \mathcal{L}}{\partial \mathbf{u}}, \delta \mathbf{u} \right \rangle &= \left \langle \frac{\partial \mathcal{J}}{\partial \mathbf{u}}, \delta \mathbf{u} \right \rangle - \left \langle \boldsymbol{\lambda}, \frac{\partial \mathcal{F}}{\partial \mathbf{u}} \delta \mathbf{u} \right \rangle = \bigg \langle \frac{\partial \mathcal{J}}{\partial \mathbf{u}} - \frac{\partial \mathcal{F}}{\partial \mathbf{u}}^\dagger \boldsymbol{\lambda}, \delta \mathbf{u} \bigg \rangle = 0 \quad \forall \, \delta \mathbf{u}, \label{eq:KKT1} \\
\left \langle \frac{\partial \mathcal{L}}{\partial \mathbf{c}}, \delta \mathbf{c} \right \rangle &= \left \langle \frac{\partial \mathcal{J}}{\partial \mathbf{c}}, \delta \mathbf{c} \right \rangle - \left \langle \boldsymbol{\lambda}, \frac{\partial \mathcal{F}}{\partial \mathbf{c}} \delta \mathbf{c} \right \rangle = \bigg \langle \frac{\partial \mathcal{J}}{\partial \mathbf{c}} - \frac{\partial \mathcal{F}}{\partial \mathbf{c}}^\dagger \boldsymbol{\lambda}, \delta \mathbf{c} \bigg \rangle = 0 \quad \forall \, \delta \mathbf{c}, \label{eq:KKT2} \\
\left \langle \frac{\partial \mathcal{L}}{\partial \boldsymbol{\lambda}}, \delta \boldsymbol{\lambda} \right \rangle &= - \left \langle \delta \boldsymbol{\lambda}, \mathcal{F} \right \rangle = 0 \quad \forall \, \delta \boldsymbol{\lambda}, \label{eq:KKT3}
\end{align} \label{eq:StationarityConditions2}%
\end{subequations}
where we have defined the adjoint $\mathcal{A}^\dagger$ of a linear operator $\mathcal{A}$ as
\begin{equation}
\langle \mathbf{a}, \mathcal{A} \mathbf{b} \rangle = \langle \mathcal{A}^\dagger \mathbf{a}, \mathbf{b} \rangle \quad \forall \, \mathbf{a}, \mathbf{b},
\end{equation}
where $\mathbf{a}$ satisfies the boundary conditions carried by the operator $\mathcal{A}$. The process of finding the adjoint operator $\mathcal{A}^\dagger$ involves integration by part and yields terminal and boundary conditions for the adjoint field $\mathbf{b}$. Thus, satisfying \eqref{eq:KKT1} for given $\mathbf{u}$ and $\mathbf{c}$ gives the adjoint equation
\begin{equation}
\frac{\partial \mathcal{J}(\mathbf{u},\mathbf{c})}{\partial \mathbf{u}} - \frac{\partial \mathcal{F}[\mathbf{u},\mathbf{c}]}{\partial \mathbf{u}}^\dagger \boldsymbol{\lambda} = 0,
\label{eq:AdjointPDE}
\end{equation}
for the adjoint field $\boldsymbol{\lambda}$, with associated terminal and boundary conditions. The third stationary condition \eqref{eq:KKT3} simply enforces the governing equation \eqref{eq:PDE} for $\mathbf{u}$ given $\mathbf{c}$, that is,
\begin{equation}
\mathcal{F}[\mathbf{u},\mathbf{c}] = 0,
\label{eq:ForwardPDE}
\end{equation}
with associated initial and boundary conditions. When \eqref{eq:KKT1} and \eqref{eq:KKT3} are satisfied, we have $\mathcal{J} = \mathcal{L}$, and \eqref{eq:KKT2} therefore gives the total gradient of the cost objective with respect to the control $\mathbf{c}$,
\begin{equation}
\frac{\mathrm{d}\mathcal{J}(\mathbf{u},\mathbf{c})}{\mathrm{d}\mathbf{c}} = \frac{\partial \mathcal{L}(\mathbf{u},\mathbf{c})}{\partial \mathbf{c}} = \frac{\partial \mathcal{J}(\mathbf{u},\mathbf{c})}{\partial \mathbf{c}} - \frac{\partial \mathcal{F}[\mathbf{u},\mathbf{c}]}{\partial \mathbf{c}}^\dagger \boldsymbol{\lambda}.
\label{eq:CostGradient}
\end{equation}
For the optimal solution, $\mathrm{d}\mathcal{J}(\mathbf{u^*},\mathbf{c^*})/\mathrm{d}\mathbf{c} =0$ holds.

There exists various adjoint-based algorithms for obtaining the optimal solution $\mathbf{u}^*, \mathbf{c}^*, \boldsymbol{\lambda}^*$ to the PDE-constrained optimization problem \eqref{eq:OptimalControlProblem}. These algorithms solve the same set of equations, namely the direct (forward) PDE \eqref{eq:ForwardPDE} and adjoint PDE \eqref{eq:AdjointPDE}, to determine the sensitivity of the cost function to the design parameters, given by \eqref{eq:CostGradient}. The difference is, however, in the manner by which the optimal solution is obtained by each algorithm.
In this work, we use the direct-adjoint-looping (DAL) algorithm \cite{foures2014optimal, nabi2017adjoint,nabi2019nonlinear}, which proceeds as follows. At each iteration $k$, one first solves the forward PDE \eqref{eq:ForwardPDE} for $\mathbf{u}^k$, given the current control $\mathbf{c}^k$. With $\mathbf{u}^k$ and $\mathbf{c}^k$ in hand, one then solves the adjoint PDE \eqref{eq:AdjointPDE} for $\boldsymbol{\lambda}^k$ in backward time since the adjoint PDE contains a terminal condition instead of an initial condition. Finally, one computes the gradient of the cost objective using \eqref{eq:CostGradient}, which is then used to update the control as
\begin{equation}
\mathbf{c}^{k+1} = \mathbf{c}^k - \beta \frac{\mathrm{d}\mathcal{J}(\mathbf{u}^k,\mathbf{c}^k)}{\mathrm{d}\mathbf{c}},
\label{eq:GradientUpdateDAL}
\end{equation}
with $\beta$ a learning rate that we will keep fixed. It should be noted that the convergence rate can be increased by employing more sophisticated update formulas such as quasi-Newton methods \cite{nocedal2006numerical}. In our case, every gradient update only requires two PDE solutions, one for the forward PDE and one for the adjoint PDE. We end the iterations once the cost objective has stopped decreasing, or alternatively once the gradient \eqref{eq:CostGradient} becomes small enough.

\subsection{Guidelines for training and evaluating the PINN optimal solution}
\label{sec:Guidelines}

The PINN methodology presented in Section \ref{sec:PINNs} involves training neural networks that approximate the state variable $\mathbf{u}$ and the control variable $\mathbf{c}$ by minimizing \textcolor{black}{the} augmented loss function \eqref{eq:LossControl}, \textcolor{black}{which can be written as
\begin{equation}
\mathcal{L}(\boldsymbol{\theta}_\mathbf{u},\boldsymbol{\theta}_\mathbf{c}) = \mathcal{L}_{\mathcal{F}/\mathcal{B}/\mathcal{I}}(\boldsymbol{\theta}_\mathbf{u},\boldsymbol{\theta}_\mathbf{c}) + w_\mathcal{J} \mathcal{L}_\mathcal{J}(\boldsymbol{\theta}_\mathbf{u},\boldsymbol{\theta}_\mathbf{c}),
\label{eq:LossControlSimp}
\end{equation}
where $\mathcal{L}_{\mathcal{F}/\mathcal{B}/\mathcal{I}}(\boldsymbol{\theta}_\mathbf{u},\boldsymbol{\theta}_\mathbf{c})$ includes the first three terms in \eqref{eq:LossControl} relating to the PDE residual, boundary conditions and initial condition. Thus, \eqref{eq:LossControlSimp} is a multi-objective loss
}
that seeks to simultaneously enforce the governing PDE and decrease the cost objective. 
However, incorporating various terms in the loss function makes the neural networks harder to train \citep{krishnapriyan2021characterizing}, and the quality of the optimal solution is not guaranteed since the cost objective might be minimized at the expense of satisfying the PDE \textcolor{black}{or vice-versa}. Issues arising from a composite loss function are not new to PINNs -- even forward problems are characterized by loss functions containing various components accounting for the PDE residual and initial/boundary conditions; see \eqref{eq:Loss}. Yet, optimal control problems pose an additional challenge due to the possible scale difference between the PDE loss \textcolor{black}{$\mathcal{L}_{\mathcal{F}/\mathcal{B}/\mathcal{I}}$} which should ideally vanish, and the cost objective \textcolor{black}{loss $\mathcal{L}_\mathcal{J}$} which might remain finite in the true optimal solution.
Several recent studies try to address the issue of balancing different objectives when training PINNs, using adaptive weighting strategies \cite{wang2020understanding,van2020optimally,maddu2022inverse,bischof2021multi}, augmented Lagrangian methods \cite{lu2021physics,basir2021physics}, \textcolor{black}{or bi-level approaches that decouple the different objectives \cite{wang2021fast,hao2022bi}}. 

\begin{figure}
\centering
\includegraphics[width=\textwidth]{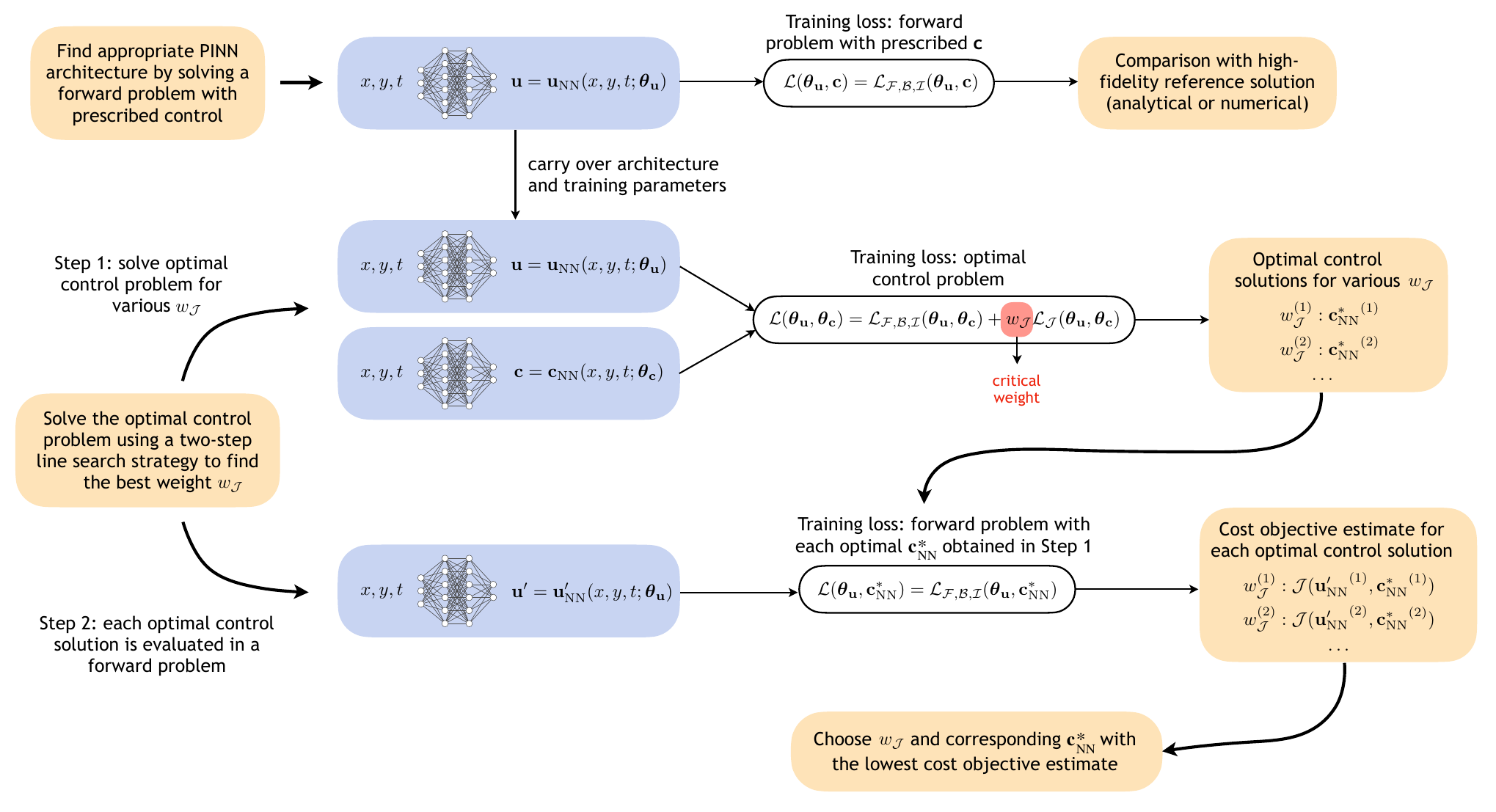}
\caption{Guidelines for obtaining a good optimal control solution using the PINN framework.}
\label{fig:Guidelines}
\end{figure}

In the present paper, however, our goal is to evaluate the feasibility and performance of the original PINN framework in solving optimal control problems. We therefore leave aside \textcolor{black}{the aforementioned} recent advances, and instead \textcolor{black}{propose simple yet effective guidelines for obtaining an optimal solution that minimizes the cost objective while satisfying the PDE constraint. These guidelines are illustrated in Figure \ref{fig:Guidelines} and consist of two parts:}
{\color{black}
\begin{enumerate}
\item First, we ensure that the network architecture employed for $\mathbf{u}_\mathrm{NN}$ and training procedure are appropriate for the type of PDE constraining the optimal control problem of interest. This can be accomplished by first computing a PINN solution to a forward problem based on the same PDE, ideally with a known solution, using the forward loss \eqref{eq:Loss} (which is equal to $\mathcal{L}_{\mathcal{F}/\mathcal{B}/\mathcal{I}}$ for a prescribed control $\mathbf{c}$). During this process, one finds a network architecture, distribution of residual points, training hyperparameters (number of epochs, batch size, etc), and weights $w_b$, $w_0$ that are tailored for the PDE at hand and can be then carried over to the solution of the optimal control problem. We emphasize that the PINN solution of forward problems governed by PDEs has been extensively treated in the literature \cite{cai2021physics,cuomo2022scientific}.

\item For the PINN solution of the optimal control problem, the only remaining parameters to choose are the network architecture for $\mathbf{c}_\mathrm{NN}$ and the weight $w_\mathcal{J}$ in the loss \eqref{eq:LossControlSimp}. In particular, the choice of $w_\mathcal{J}$ is critical since it controls a trade-off between the two potentially conflicting objectives $\mathcal{L}_{\mathcal{F}/\mathcal{B}/\mathcal{I}}$ and $\mathcal{L}_\mathcal{J}$: large $w_\mathcal{J}$ means that $\mathcal{L}_\mathcal{J}$ is minimized at the expense of $\mathcal{L}_{\mathcal{F}/\mathcal{B}/\mathcal{I}}$, while small $w_\mathcal{J}$ leads to the opposite behavior. To find the optimal $w_\mathcal{J}$ yielding a solution that simultaneously minimizes the cost objective and satisfies the PDE constraint, we propose a two-step line search strategy.

\textit{Step 1.} Following the PINN methodology presented in Section \ref{sec:PINNs}, we compute a solution $\mathbf{c}_\mathrm{NN}^*$ to the optimal control problem by creating a first neural network $\mathbf{u}_\mathrm{NN}$ for the state variable with the architecture obtained previously, a second neural network $\mathbf{c}_\mathrm{NN}$ for the control variable, and by minimizing the loss \eqref{eq:LossControlSimp}. We repeat this process for a range of values of $w_\mathcal{J}$.

\textit{Step 2.} We evaluate the performance of the optimal control $\mathbf{c}_\mathrm{NN}^*$ obtained for each $w_\mathcal{J}$ by solving the forward problem corresponding to $\mathbf{c}_\mathrm{NN}^*$ with another PINN. Specifically, we train a separate neural network $\mathbf{u}_\mathrm{NN}'$ to minimize the loss \eqref{eq:Loss} with fixed $\mathbf{c} = \mathbf{c}_\mathrm{NN}^*$ (which corresponds to $\mathcal{L}_{\mathcal{F}/\mathcal{B}/\mathcal{I}}$ with fixed control). We then use the solution to this forward problem to compute an approximate cost objective $\mathcal{J}(\mathbf{u}_\mathrm{NN}',\mathbf{c}_\mathrm{NN}^*)$. Finally, we select as our final solution the optimal control $\mathbf{c}_\mathrm{NN}^*$ given by the value of $w_\mathcal{J}$ that yields the smallest approximate $\mathcal{J}$.
\end{enumerate}
}

\section{Results}
\label{sec:Results}

We apply the PINN and adjoint-based approaches for solving optimal control problems to four PDEs modeling a range of physical systems -- the Laplace, Burgers, Kuramoto-Sivashinsky, and Navier-Stokes equations. In each case, we follow \textcolor{black}{the guidelines presented in Section \ref{sec:Guidelines} to obtain the PINN solution. That is, we first validate our PINN architecture by solving a forward problem without (or with given) control and comparing the PINN solution with a reference solution (which is either obtained analytically or from a high-fidelity numerical code). We then solve the optimal control problem using the two-step line search strategy described above.}

\textcolor{black}{We compare rigorously the performance of the optimal control solutions $\mathbf{c}^*$ obtained from the PINN and adjoint-based DAL frameworks, which is quantified by the value of the cost objective $\mathcal{J}$ corresponding to $\mathbf{c}^*$. To compute an accurate estimate of $\mathcal{J}$, we use a high-fidelity numerical code (based on finite-element or spectral methods) to compute the solution $\mathbf{u}_\mathrm{HF}$ to the forward problem with fixed control equal to $\mathbf{c}^*$ coming from the PINN or DAL solution, yielding an accurate cost objective estimate $\mathcal{J}(\mathbf{u}_\mathrm{HF},\mathbf{c}^*)$. This ensures a rigorous evaluation and fair comparison of the PINN and DAL optimal control solutions $\mathbf{c}^*$.} 

\subsection{Laplace equation}
\label{sec:Laplace}

Let us first consider the Laplace equation
\begin{equation}
\frac{\partial^2 u}{\partial x^2} + \frac{\partial^2 u}{\partial y^2} = 0,
\label{eq:Laplace}
\end{equation}
where the potential $u(x,y)$ is defined in a square domain $(x,y) \in [0,1] \times [0,1]$. The boundary conditions will be specified below. The linearity of the Laplace equation implies that any derivative optimal control problem defined by a quadratic cost objective will be convex \cite{troltzsch2010optimal}, which makes it an ideal setting for a first comparison between the PINN and DAL frameworks.

\subsubsection{Forward problem}

We first solve the forward problem defined by the boundary conditions
\begin{equation}
u(x,1) = \sin \pi x, \quad u(x,0) = u(0,y) = u(1,y) = 0.
\label{eq:LaplaceForwardBC}
\end{equation}
For these boundary conditions, the Laplace equation admits the analytical solution
\begin{equation}
u_a(x,y) = \cos(\pi x) \sinh(\pi y),
\label{eq:LaplaceForwardExact}
\end{equation}
which we will use as a benchmark to evaluate the PINN solution.

To solve the problem in the PINN framework, we represent $u(x,y)$ with a neural network containing 4 hidden layers of 50 neurons each, which we train using the loss \eqref{eq:Loss} according to the gradient update formula \eqref{eq:GradientUpdate}. To evaluate the loss, we sample $N_r = 10000$ residual training points $(x_i,t_i) \in [0,1] \times [0,1]$ using a Latin hypercube sampling (LHS) strategy, and we select $N_b = 160$ equally-spaced boundary training points $(x_i,t_i)$ on the boundary of the domain. The entire set of $N_r$ residual points is randomly separated into 10 minibatches of $N_r/10 = 1000$ points each. During each gradient update \eqref{eq:GradientUpdate}, the residual loss component is calculated using one minibatch of residual points. One epoch of training is defined as a complete pass through the entire set of $N_r$ residual points, that is, through all 10 minibatches. At the beginning of each epoch, the $N_r$ residual points are shuffled and separated into 10 new minibatches. We set uniform scalar weights $w_r = w_b = 1$. Finally, we choose an initial learning rate of $\alpha = 10^{-3}$ and decrease it by a factor 10 after 3000 epochs, for a total of 6000 training epochs.

\begin{figure}
\centering
\includegraphics[width=\textwidth]{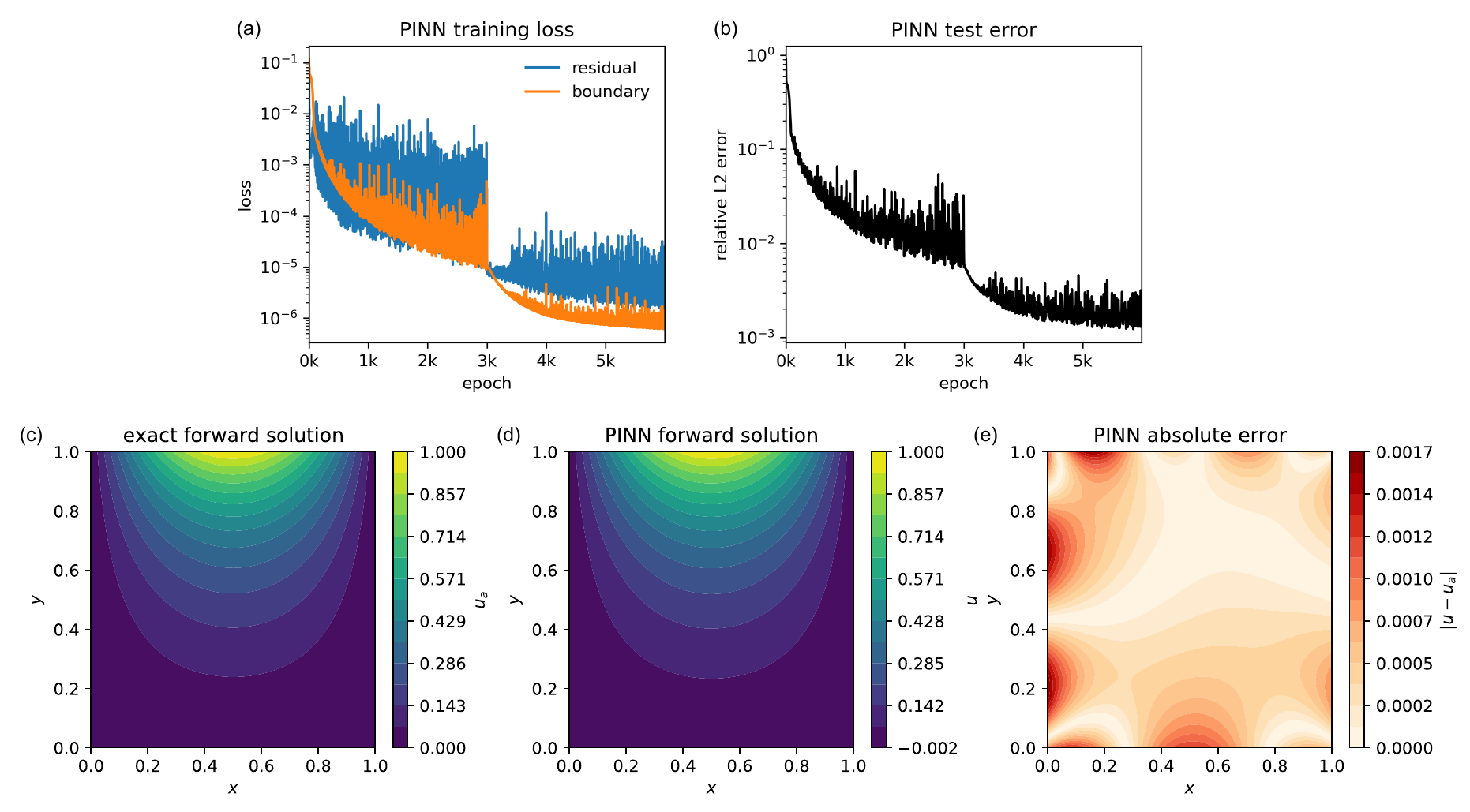}
\caption{Forward solution of the Laplace equation \eqref{eq:Laplace}. (a,b) Convergence of the loss and relative $L_2$ test error during training of the PINN solution. (c) Analytical solution $u_a$, (d) trained PINN solution $u$ and (e) its local absolute error $|u-u_a|$.}
\label{fig:LaplaceForward}
\end{figure}

The different loss components during the training process are shown in Figure \ref{fig:LaplaceForward}(a), and the relative $L_2$ error of the PINN solution versus its analytical counterpart, $||u-u_a||_2/||u_a||_2$, is displayed in Figure \ref{fig:LaplaceForward}(b).
The $L_2$ error is equivalent to a test error since it is computed on a $100 \times 100$ cartesian grid independent from the residual training points. Figures \ref{fig:LaplaceForward}(a,b) reveal that decreasing loss components are directly correlated with a decreasing test error, and can thus be monitored during training to gauge the accuracy of the PINN solution. The analytical solution $u_a$ and the trained PINN solution $u$ are displayed in Figures \ref{fig:LaplaceForward}(c,d), while the point-wise absolute difference between the two is reported in Figure \ref{fig:LaplaceForward}(e). The excellent accuracy of the PINN forward solution validates the choice of neural network parameters and training points.

\subsubsection{Optimal control problem}

Next, we use the same PINN architecture to solve an optimal control problem constrained on the Laplace equation. Consider now the boundary conditions
\begin{equation}
u(x,0) = \sin \pi x, \quad u(x,1) = f(x), \quad u(0,y) = u(1,y), \quad \frac{\partial u}{\partial x}(0,y) = \frac{\partial u}{\partial x}(1,y)
\label{eq:LaplaceControlBC}
\end{equation}
where $f(x)$ is a control potential applied to the top wall, \textcolor{black}{and the boundary conditions are now periodic along the $x$ direction}. We then seek to solve the convex optimal control problem defined as
\begin{equation}
f^* = \arg \min_{f} \mathcal{J}(u) \quad \text{subject to } \eqref{eq:Laplace} \text{ and } \eqref{eq:LaplaceControlBC},
\label{eq:OptimalControlLaplace}
\end{equation}
where the cost objective is
\begin{equation}
\mathcal{J}(u) = \int_0^1 \left| \frac{\partial u}{\partial y}(x,1) - q_d(x) \right|^2 dx, \qquad q_d(x) = \cos(\pi x).
\label{eq:CostLaplace}
\end{equation}
In other words, we want to find the optimal potential $f^*(x)$ at the top wall that produces the desired flux $q_d(x)$. This problem has the analytical optimal solution
\begin{subequations}
\begin{align}
f_a^*(x) &= \mathrm{sech} (2 \pi) \sin (2 \pi x) + \frac{1}{2 \pi} \mathrm{tanh} (2 \pi) \cos (2 \pi x), \\
u_a^*(x,y) &= \frac{1}{2} \mathrm{sech} (2 \pi) \sin (2 \pi x) (e^{2 \pi (y-1)} + e^{2 \pi (1-y)}) + \frac{1}{4 \pi} \mathrm{sech} (2 \pi) \cos (2 \pi x) (e^{2 \pi y} - e^{-2 \pi y}),
\end{align} \label{eq:LaplaceControlAnalytical}
\end{subequations}
which we will use to evaluate the accuracy of the PINN and DAL optimal solutions. 

To solve the optimal control problem with PINNs, we define a second neural network for $f(x)$ consisting of 3 hidden layers of 30 neurons each. Following the guidelines presented in Section \ref{sec:Guidelines}, we keep the same architecture as before (albeit with newly initialized parameters) for the network representing $u(x,y)$. We then train simultaneously both neural networks using the loss \eqref{eq:LossControl} and the gradient update formula \eqref{eq:GradientUpdateControl}. We use the same training points as before, and evaluate the integral in the cost objective \eqref{eq:CostLaplace} using the midpoint rule at $N_\mathcal{J} = 41$ equally-spaced training points $x_i \in [0,1]$. We use uniform scalar weights $w_r = w_b = 1$, we choose an initial learning rate of $\alpha = 10^{-3}$ and decrease it by a factor 10 after 5000 epochs, for a total of 10000 training epochs. \textcolor{black}{We repeat this procedure for 11 values of $w_\mathcal{J}$ between $10^{-3}$ and $10^7$.}

To find the DAL optimal solution, we implement the DAL iterative procedure in the finite-volume solver OpenFOAM. The adjoint Laplace equation and the gradient of the cost objective are given in \cite{hazra2004simultaneous}, in which the same problem is considered. The DAL optimal solution is obtained by iteratively solving the Laplace equation and its adjoint on a $40 \times 40$ grid, updating the control $f(x)$ at each iteration with the gradient descent formula \eqref{eq:GradientUpdateDAL}. We start the iterations with a zero initial guess for $f(x)$ and employ a learning rate $\beta = 1$.

\textcolor{black}{Beginning with Step 1 of the line search strategy presented in Section \ref{sec:Guidelines}, we visualize in Figure \ref{fig:LaplaceControl}(a) the two components $\mathcal{L}_{\mathcal{F}/\mathcal{B}}$ and $\mathcal{L}_{\mathcal{J}}$ of the loss \eqref{eq:LossControlSimp} obtained at the end of training of the PINN optimal control solution, for each considered value of $w_\mathcal{J}$. As expected, $\mathcal{L}_{\mathcal{F}/\mathcal{B}}$ and $\mathcal{L}_{\mathcal{J}}$ vary in opposite directions as $w_\mathcal{J}$ increases.
\begin{figure}[t]
\centering
\includegraphics[width=\textwidth]{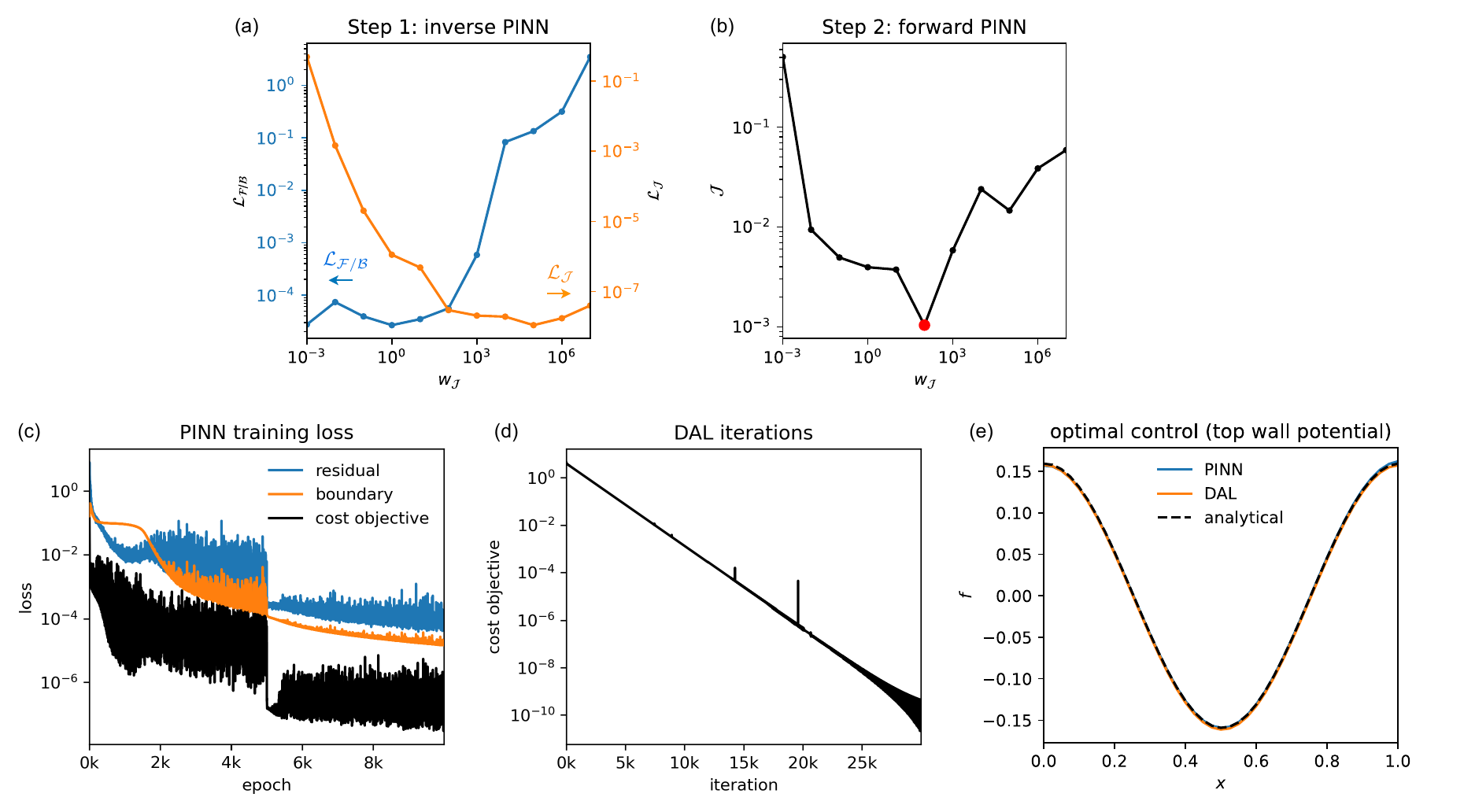}
\caption{Optimal solution of the Laplace control problem \eqref{eq:OptimalControlLaplace}. \textcolor{black}{(a) Components of the loss \eqref{eq:LossControlSimp} obtained at the training of the PINN optimal control solution versus weight $w_\mathcal{J}$ (step 1 of the line search strategy). (b) Cost objective estimate obtained by a separate PINN solution of the forward problem with fixed control from the PINN optimal solution versus $w_\mathcal{J}$ (step 2 of the line search strategy). The best optimal control, obtained with $w_\mathcal{J} = 100$, is shown by the red dot.} (c) Convergence of the loss during training of the PINN \textcolor{black}{optimal control solution (for $w_\mathcal{J} = 100$)}. (d) Convergence of the cost objective during DAL iterations. (e) Optimal top wall potential $f^*$ obtained from the PINN \textcolor{black}{(for $w_\mathcal{J} = 100$)} and DAL frameworks, compared with the top wall potential $f_a^*$ of the analytical solution.}
\label{fig:LaplaceControl}
\end{figure}
Moving to Step 2 of the line search strategy, we report in Figure \ref{fig:LaplaceControl}(b) the cost objective value $\mathcal{J}$ obtained by a separate PINN solution of the corresponding forward problem, with fixed control from the optimal solution obtained in Step 1, for each value of $w_\mathcal{J}$. We observe that the control solution obtained with $w_\mathcal{J} = 100$, shown by the red dot, yields the lowest $\mathcal{J}$. It is therefore this optimal control solution that we analyze hereafter and compare with the DAL solution.}

The various loss components during training of the PINN optimal control solution \textcolor{black}{(for $w_\mathcal{J} = 100$)} are shown in Figure \ref{fig:LaplaceControl}(c), while the convergence of the cost objective during the DAL iterations is shown in Figure \ref{fig:LaplaceControl}(d). Figure \ref{fig:LaplaceControl}(e) shows the optimal top wall potentials $f^*$ obtained from the PINN and DAL frameworks, together with their analytical counterpart $f_a^*$. The excellent agreement between the three potentials demonstrates that the PINN and DAL frameworks both found the global solution to this convex optimal control problem.

\subsection{Burgers equation}
\label{sec:Burgers}

We then consider the one-dimensional Burgers equation, a prototypical nonlinear hyperbolic PDE that takes the form
\begin{equation}
\frac{\partial u}{\partial t} + u \frac{\partial u}{\partial x} = \nu \frac{\partial^2 u}{\partial x^2}, 
\label{eq:Burgers}
\end{equation}
where $u(x,t)$ is the velocity at position $x \in [0,L]$ and time $t \in [0,T]$. We assign periodic boundary conditions and select the viscosity $\nu = 0.01$; the initial condition $u(x,0) = u_0(x)$ will be specified later. We take $L = 4$, for which the Burgers equation possesses the analytical solution
\begin{equation}
u_a(x,t) = \frac{2 \nu \pi e^{-\pi^2 \nu (t-5)} \sin(\pi x)}{2 + e^{-\pi^2 \nu (t-5)} \cos(\pi x)}.
\label{eq:BurgersAnalytical}
\end{equation}

\subsubsection{Forward problem}

For the forward problem, we compute the evolution of the system up to the time horizon $T = 5$ given the initial condition $u_0(x) = u_a(x,0)$. In this way, the accuracy of the PINN solution can be evaluated against the analytical solution \eqref{eq:BurgersAnalytical}. 

The PINN solution is obtained by representing $u(x,t)$ with a network containing 4 hidden layers of 50 neurons each, which we train using the loss \eqref{eq:Loss}. To evaluate the loss, we sample $N_r = 20000$ residual training points $(x_i,t_i) \in [0,L] \times [0,T]$ using LHS. We select $N_b = 82$ equally-spaced boundary training points $(x_i,t_i) \in \{0,L\} \times [0,T]$, and $N_0 = 41$ equally-spaced initial training points $(x_i,t_i) \in [0,L] \times \{0\}$. As before, the set of $N_r$ residual points is randomly separated into 10 minibatches of $N_r/10 = 2000$ points, and we set uniform scalar weights $w_r = w_b = w_0 = 1$. Finally, we choose an initial learning rate of $\alpha = 10^{-3}$ and decrease it by a factor 10 after 5k epochs of training, for a total of 10k epochs.

\begin{figure}
\centering
\includegraphics[width=\textwidth]{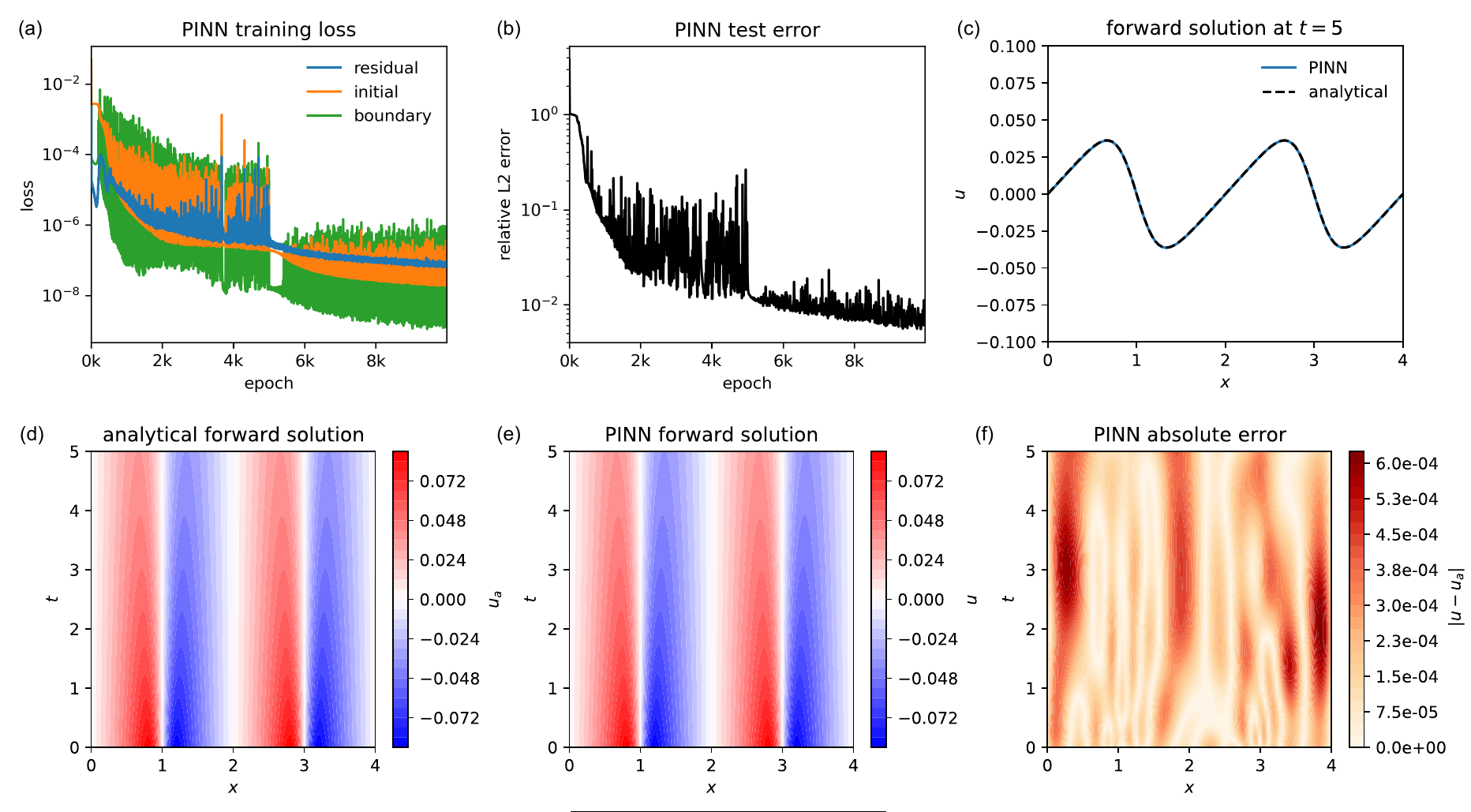}
\caption{Forward solution of the Burgers equation \eqref{eq:Burgers}. (a,b) Convergence of the loss and relative $L_2$ test error during training of the PINN solution. (c) Snapshots at final time of the trained PINN and analytical solutions. (d) Analytical solution $u_a$, (e) trained PINN solution $u$ and (f) its local absolute error $|u-u_a|$.}
\label{fig:BurgersForward}
\end{figure}

The different loss components during the training process are shown in Figure \ref{fig:BurgersForward}(a), and the relative $L_2$ error of the PINN solution versus its analytical counterpart, $||u-u_a||_2/||u_a||_2$, is displayed in Figure \ref{fig:BurgersForward}(b).
As was observed in the case of the Laplace equation, decreasing loss components are directly correlated with a decreasing solution error. Snapshots of the trained PINN and analytical solutions at final time $t = 5$, shown in Figure \ref{fig:BurgersForward}(c), are in excellent agreement. This is confirmed by the contour plots of the analytical and trained PINN solutions displayed in Figures \ref{fig:BurgersForward}(d,e), and the absolute error between them shown in Figure \ref{fig:BurgersForward}(f). The good accuracy of the PINN solution validates the choice of neural network parameters and training points.

\subsubsection{Optimal control problem}

Next, we turn our attention to an optimal control problem, which we define as
\begin{equation}
u_0^* = \arg \min_{u_0} \mathcal{J}(u) \quad \text{subject to } \eqref{eq:Burgers},
\label{eq:OptimalControlBurgers}
\end{equation}
where the objective cost is
\begin{equation}
\mathcal{J}(u) = \frac{1}{2} \int_0^L |u(x,T) - u_a(x,T)|^2 dx.
\label{eq:CostBurgers}
\end{equation}
In other words, we seek the optimal initial condition $u_0^*(x)$ that results in the same final state as the analytical solution \eqref{eq:BurgersAnalytical}. Because of the nonlinearity of the Burgers equation, this optimal control problem is no longer convex as was the case for the Laplace equation, and it may therefore possess multiple local minima in addition to the global minimum given by $u_a(x,0)$.

The PINN optimal solution is obtained by defining a second neural network for $u_0(x)$, this time consisting of 3 hidden layers of 30 neurons each, in addition to the network for $u(x,t)$ for which we keep the same architecture as before. We then train simultaneously both neural networks using the loss \eqref{eq:LossControl}, starting from a new initialization of the parameters. We use the same training points as before, and evaluate the integral in the cost objective \eqref{eq:CostBurgers} using the midpoint rule at $N_\mathcal{J} = 41$ equally-spaced cost training points $x_i \in [0,L]$. We use uniform scalar weights $w_r = w_b = w_0 = w_\mathcal{J} = 1$, choose an initial learning rate of $\alpha = 10^{-3}$ and decrease it by a factor 10 after 20k and 25k epochs of training, for a total of 30k epochs. \textcolor{black}{We repeat this procedure for 10 values of $w_\mathcal{J}$ between $10^{-3}$ and $10^6$.}

The DAL optimal solution is obtained by iteratively solving the Burgers and adjoint Burgers equations, updating the control $u_0(x)$ at each iteration with the gradient descent formula \eqref{eq:GradientUpdateDAL}. The adjoint Burgers equation and the gradient of the cost objective are given in \ref{sec:AdjointBurgersEquation}. A spectral solver with 256 Fourier modes and semi-implicit Euler scheme with $dt = 10^{-3}$ is used to solve the forward and adjoint Burgers equations at each iteration. We start the iterations with a zero initial guess for the control $u_0(x)$ and employ a learning rate $\beta = 1$.

\textcolor{black}{Beginning with Step 1 of the line search strategy presented in Section \ref{sec:Guidelines}, we visualize in Figure \ref{fig:BurgersControl}(a) the two components $\mathcal{L}_{\mathcal{F}/\mathcal{B}/\mathcal{I}}$ and $\mathcal{L}_{\mathcal{J}}$ of the loss \eqref{eq:LossControlSimp} obtained at the end of training of the PINN optimal control solution, for each considered value of $w_\mathcal{J}$.
\begin{figure}
\centering
\includegraphics[width=\textwidth]{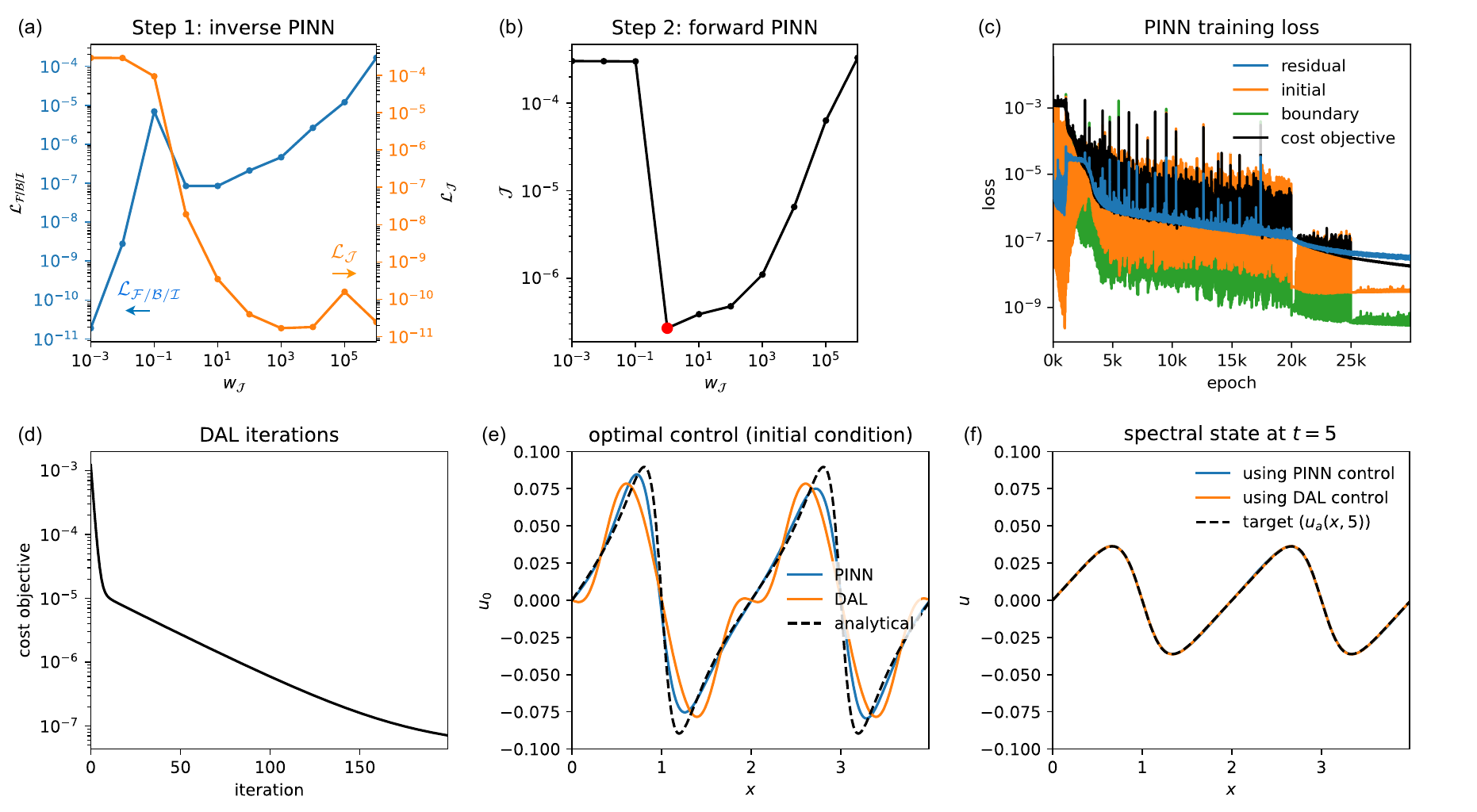}
\caption{Optimal solution of the Burgers control problem \eqref{eq:OptimalControlBurgers}. \textcolor{black}{(a) Components of the loss \eqref{eq:LossControlSimp} obtained at the end of training of the PINN optimal control solution versus weight $w_\mathcal{J}$ (step 1 of the line search strategy). (b) Cost objective estimate obtained by a separate PINN solution of the forward problem with fixed control from the PINN optimal solution versus $w_\mathcal{J}$ (step 2 of the line search strategy). The best optimal control, obtained with $w_\mathcal{J} = 1$, is shown by the red dot.} (c) Convergence of the loss during training of the PINN \textcolor{black}{optimal control solution (for $w_\mathcal{J} = 1$)}. (d) Convergence of the cost objective during DAL iterations. (e) Optimal initial condition $u_0^*$ obtained using the PINN \textcolor{black}{(for $w_\mathcal{J} = 1$)} and DAL frameworks, compared with the initial state of the analytical solution. (f) Snapshots at final time of two spectral solutions calculated using the optimal initial conditions $u_0^*$ from the PINN Convergence of the loss during training of the PINN optimal control solution \textcolor{black}{(for $w_\mathcal{J} = 1$)} and DAL frameworks, compared with the final state of the analytical solution.}
\label{fig:BurgersControl}
\end{figure}
Moving to Step 2 of the line search strategy, we report in Figure \ref{fig:BurgersControl}(b) the cost objective value $\mathcal{J}$ obtained by a separate PINN solution of the corresponding forward problem, with fixed control from the optimal solution obtained in Step 1, for each value of $w_\mathcal{J}$. We observe that the control solution obtained with $w_\mathcal{J} = 1$, shown by the red dot, yields the lowest $\mathcal{J}$. It is therefore this optimal control solution that we analyze hereafter and compare with the DAL solution.}

The various loss components during training of the PINN optimal solution \textcolor{black}{(for $w_\mathcal{J} = 1$)} are displayed in Figure \ref{fig:BurgersControl}(c), while the convergence of the cost objective during the DAL iterations is shown in Figures \ref{fig:BurgersControl}(d). Figure \ref{fig:BurgersControl}(e) compares the optimal initial conditions $u_0^*$ obtained using the PINN and DAL frameworks, together with the initial state $u_a(x,0)$ of the analytical solution \eqref{eq:BurgersAnalytical}. In order to evaluate the quality of the PINN and DAL optimal solutions, we use a spectral solver (the same used in the DAL iterations) to compute the final state $u(x,5)$ corresponding to the optimal initial conditions $u_0^*$ from the PINN and DAL frameworks. Figure \ref{fig:BurgersControl}(f) shows the resulting computed states $u(x,5)$ alongside the target final state $u_a(x,5)$ of the analytical solution. The corresponding cost values are $\mathcal{J} = 2.27 \cdot 10^{-7}$ and $\mathcal{J} = 7.12 \cdot 10^{-8}$ for the PINN and DAL optimal initial conditions, respectively. Although the initial conditions identified by the PINN and DAL frameworks differ substantially from $u_a(x,0)$, they both result in a final state that is almost identical to $u_a(x,5)$. This is a consequence of the nonlinearity of the Burgers equation, which allows noticeably different initial conditions to result in extremely similar final states -- in other words, the optimization landscape for this problem admits multiple local minima with small cost objective values. In this situation, it is interesting to observe that the DAL procedure finds a slightly better optimal initial condition $u_0^*$ than the PINN framework, but at the expense of its smoothness.

\subsection{Kuromoto-Sivashinsky equation}
\label{sec:KS}

The next example we consider is the one-dimensional Kuramoto-Sivashinsky (KS) equation, one of the simplest PDEs that generates chaotic behavior under certain conditions. The KS equation takes the form
\begin{equation}
\frac{\partial u}{\partial t} + u \frac{\partial u}{\partial x} + \frac{\partial^2 u}{\partial x^2} + \frac{\partial^4 u}{\partial x^4} = f(x,t), 
\label{eq:KS}
\end{equation}
where $u(x,t)$ is the velocity at position $x \in [0,L]$ and time $t \in [0,T]$, and $f(x,t)$ is a distributed control force. We use periodic boundary conditions and choose as initial condition 
\begin{equation}
u_0(x) = \cos \left( \frac{2 \pi x}{10} \right) + \mathrm{sech}  \left( \frac{x-L/2}{5} \right).
\label{eq:KSIC}
\end{equation}
Solutions to the KS equation without forcing undergo a sequence of bifurcations as $L$ increases. \textcolor{black}{The} zero state is a stable \textcolor{black}{fixed-point} solution for $L < 2\pi$, but becomes linearly unstable for $L > 2 \pi$. \textcolor{black}{For} even larger $L$, the solution becomes chaotic. Here, we choose $L = 50$ which corresponds to the chaotic regime \cite{cvitanovic2010state}, with the aim of finding a \textcolor{black}{control force} $f(x,t)$ that drives the \textcolor{black}{state $u$ towards the unstable zero fixed-point solution}.

\subsubsection{Forward problem}

For the forward problem, we solve the KS equation up to the time horizon $T = 10$ given the initial condition \eqref{eq:KSIC} and $f(x,t) = 0$. The reference solution, $u_s$, is given by a spectral discretization of the KS equation, employing 256 Fourier modes and a semi-implicit Euler scheme with time step $dt = 10^{-4}$.

The PINN solution is obtained by representing $u(x,t)$ with a network containing 5 hidden layers of 50 neurons each, which we train by minimizing the loss \eqref{eq:Loss}. To this end, we sample $N_r = 80000$ residual training points $(x_i,t_i) \in [0,L] \times [0,T]$ using LHS. We select $N_b = 82$ equally-spaced boundary training points $(x_i,t_i) \in \{0,L\} \times [0,T]$, and $N_0 = 41$ equally-spaced initial training points $(x_i,t_i) \in [0,L] \times \{0\}$. At the beginning of each epoch, the entire set of $N_r$ residual points is shuffled and divided into 20 minibatches of $N_r/20 = 4000$ points each. We set scalar weights $w_r = w_b = w_0 = 1$, choose an initial learning rate of $\alpha = 10^{-3}$ and decrease it by a factor 10 after 10k and 20k epochs, for a total of 30k epochs.

\begin{figure}
\centering
\includegraphics[width=\textwidth]{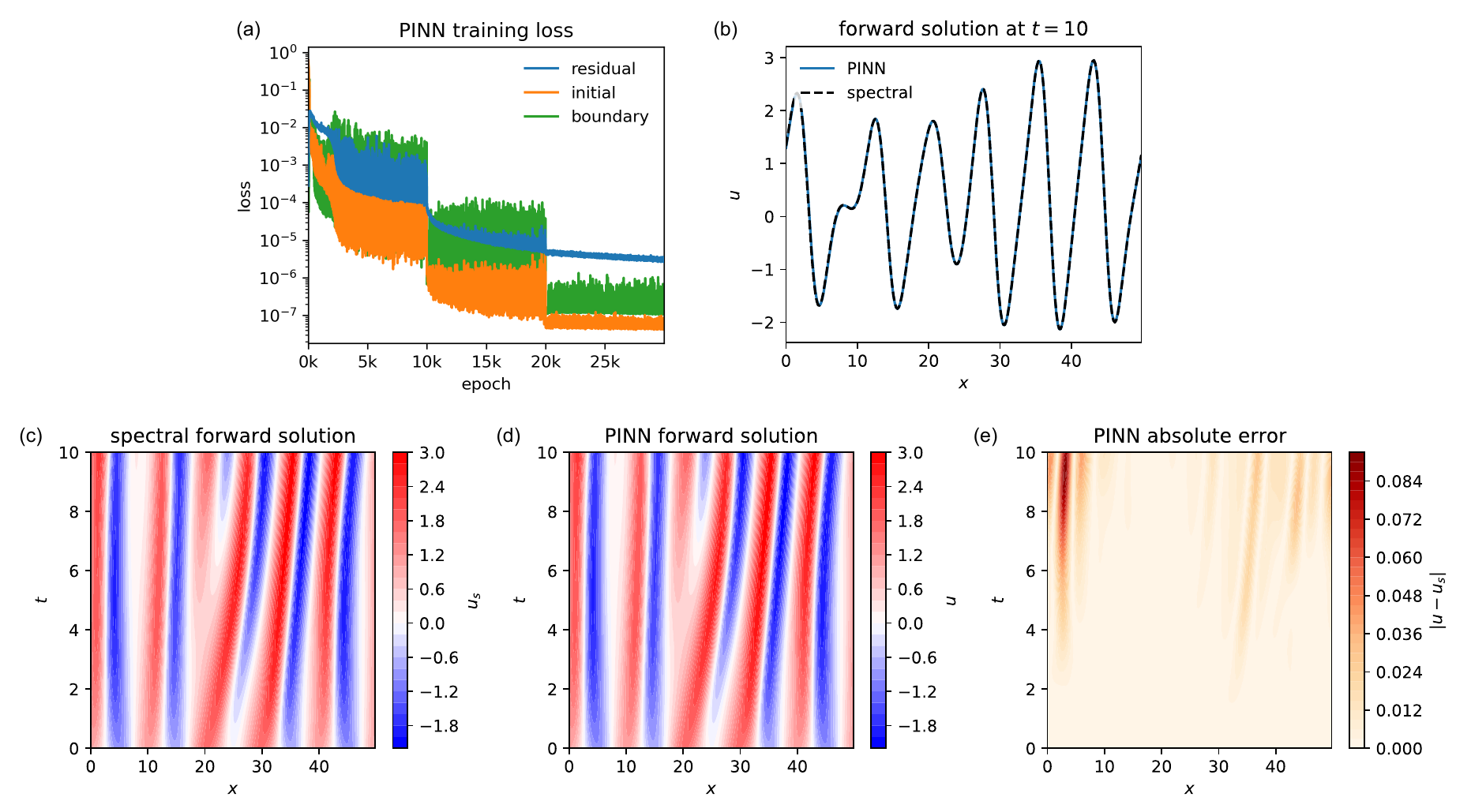}
\caption{Forward solution of the KS equation \eqref{eq:KS}. (a) Convergence of the loss during training of the PINN solution. (b) Snapshots at final time of the trained PINN solution and a reference spectral solution. (d) Reference spectral solution $u_s$, (e) trained PINN solution $u$ and (f) its local absolute error $|u-u_s|$}
\label{fig:KSForward}
\end{figure}

The different loss components during the training of the PINN forward solution are shown in Figure \ref{fig:KSForward}(a).
Snapshots of the trained PINN and reference spectral solutions at final time $t=10$ are shown in Figure \ref{fig:KSForward}(b), displaying excellent agreement with each other. Contour plots of the spectral and trained PINN solutions are displayed in Figures \ref{fig:KSForward}(c,d), and the absolute error between the two is shown in Figure \ref{fig:KSForward}(e). The low error of the PINN solution validates the choice of neural network parameters and training points.

\subsubsection{Optimal control problem}

We define the optimal control problem as
\begin{equation}
f^* = \arg \min_{f} \mathcal{J}(u,f) \quad \text{subject to } \eqref{eq:KS},
\label{eq:OptimalControlKS}
\end{equation}
where the objective cost is
\begin{equation}
\mathcal{J}(u,f) = \frac{1}{2} \int_0^T \int_0^L (|u(x,t)|^2 + \sigma |f(x,t)|^2) dx dt,
\label{eq:CostKS}
\end{equation}
\textcolor{black}{in which we choose $\sigma = 1$.} Thus, we seek the optimal \textcolor{black}{control force} $f^*(x,t)$ that drives the \textcolor{black}{system state towards the unstable zero fixed-point solution} by minimizing a quadratic cost that \textcolor{black}{penalizes} the norms of both the \textcolor{black}{state $u$} and the \textcolor{black}{control force $f$}. \textcolor{black}{Penalizing the norm of $f$ is a way to regularize the optimal control problem, without which the solution for $f^*$ would be one that instantaneously pushes the state to zero right after initial time, resulting in infinitely large $f^*$ at initial time. In this particular example, therefore, including such regularization in \eqref{eq:CostKS} is necessary to prevent the optimal control force from growing unbounded at initial time. We note that our formulation mimics the classical problem in control theory of finding a controller that drives the state of a dynamical system towards an unstable fixed point, which is usually solved by minimizing a quadratic cost functional of the same form as \eqref{eq:CostKS}.}

To solve this problem in the PINN framework, we define a second neural network for $f(x,t)$ with the same architecture as that for $u(x,t)$, which was validated in the forward problem. We then train both networks simultaneously using the loss \eqref{eq:LossControl}, starting from a new initialization of the parameters. We use the same training points as before, and evaluate the integral in the cost objective \eqref{eq:CostKS} using Monte Carlo integration with the same minibatch of residual training points used in evaluating the residual loss component. We use the scalar weights $w_r = w_b = w_0 = 1$ and $w_\mathcal{J} = 10^{-3}$, choose an initial learning rate of $10^{-3}$ and decrease it by a factor 10 after 10k and 20k epochs of training, for a total of 30k epochs. \textcolor{black}{We repeat this procedure for 10 values of $w_\mathcal{J}$ between $10^{-8}$ and $10$.}

The DAL optimal solution is obtained by iteratively solving the KS and adjoint KS equations, updating the control $f(x,t)$ at each iteration with the gradient descent formula \eqref{eq:GradientUpdateDAL}. The adjoint KS equation and the gradient of the cost objective are given in \ref{sec:AdjointKSEquation}. A spectral solver with 256 Fourier modes and semi-implicit Euler scheme with $dt = 10^{-4}$ is used to solve the forward and adjoint KS equations at each iteration. We start the iterations with a zero initial guess for the control $f(x,t)$ and employ a learning rate $\beta = 0.001$.

\textcolor{black}{Beginning with Step 1 of the line search strategy presented in Section \ref{sec:Guidelines}, we visualize in Figure \ref{fig:KSControl}(a) the two components $\mathcal{L}_{\mathcal{F}/\mathcal{B}/\mathcal{I}}$ and $\mathcal{L}_{\mathcal{J}}$ of the loss \eqref{eq:LossControlSimp} obtained at the end of training of the PINN optimal control solution, for each considered value of $w_\mathcal{J}$.
\begin{figure}
\centering
\includegraphics[width=\textwidth]{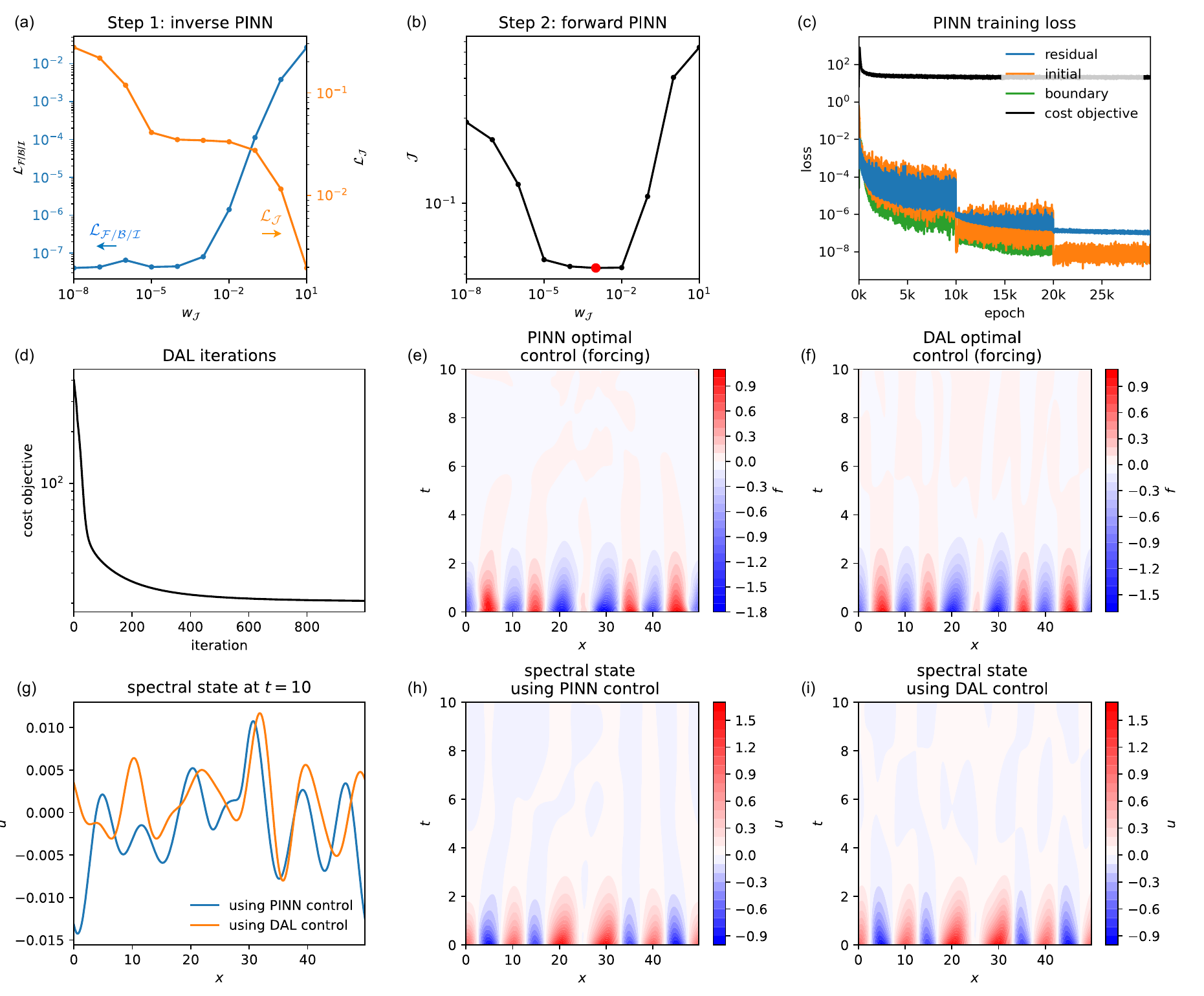}
\caption{Optimal solution of the KS control problem \eqref{eq:OptimalControlKS}. \textcolor{black}{(a) Components of the loss \eqref{eq:LossControlSimp} obtained at the end of training of the PINN optimal control solution versus weight $w_\mathcal{J}$ (step 1 of the line search strategy). (b) Cost objective estimate obtained by a separate PINN solution of the forward problem with fixed control from the PINN optimal solution versus $w_\mathcal{J}$ (step 2 of the line search strategy). The best optimal control, obtained with $w_\mathcal{J} = 10^{-3}$, is shown by the red dot.} (c) Convergence of the loss during training of the PINN \textcolor{black}{optimal control solution (for $w_\mathcal{J} = 10^{-3}$)}. (d) Convergence of the cost objective during DAL iterations. (e,f) Optimal \textcolor{black}{control forces} $f^*$ obtained from the PINN \textcolor{black}{(for $w_\mathcal{J} = 10^{-3}$)} and DAL frameworks. (g,h,i) Snapshots at final time and contour plots of two spectral solutions calculated using the optimal \textcolor{black}{control forces} $f^*$ obtained from the PINN \textcolor{black}{(for $w_\mathcal{J} = 10^{-3}$)} and DAL frameworks.}
\label{fig:KSControl}
\end{figure}
Moving to Step 2 of the line search strategy, we report in Figure \ref{fig:KSControl}(b) the cost objective value $\mathcal{J}$ obtained by a separate PINN solution of the corresponding forward problem, with fixed control from the optimal solution obtained in Step 1, for each value of $w_\mathcal{J}$. We observe that the control solution obtained with $w_\mathcal{J} = 10^{-3}$, shown by the red dot, yields the lowest $\mathcal{J}$. It is therefore this optimal control solution that we analyze hereafter and compare with the DAL solution.}

The different loss components during training of the PINN optimal solution \textcolor{black}{(for $w_\mathcal{J} = 10^{-3}$)} are displayed in Figure \ref{fig:KSControl}(c). The convergence of the cost objective during the DAL iterations is displayed in Figure \ref{fig:KSControl}(d). The optimal forcing solutions found by the PINN and DAL frameworks are shown in Figures \ref{fig:KSControl}(e,f) and look very similar. We evaluate the quality of these optimal forcings by using a spectral solver (the same used in the DAL iterations) to compute the corresponding state $u(x,t)$ up to the time horizon $T = 10$. Snapshots at final time and contour plots of the resulting states are shown in Figures \ref{fig:KSControl}(g,h,i), showing that both PINN and DAL optimal forcings manage to drive the state towards near-zero values. The corresponding cost objectives, calculated from the spectral solutions, have remarkably similar values of $\mathcal{J} = 20.58$ and $\mathcal{J} = 20.64$ for the PINN and DAL optimal forcings, respectively.

\subsection{Navier-Stokes equations}
\label{sec:NS}

\begin{figure}
\centering
\includegraphics[width=\textwidth]{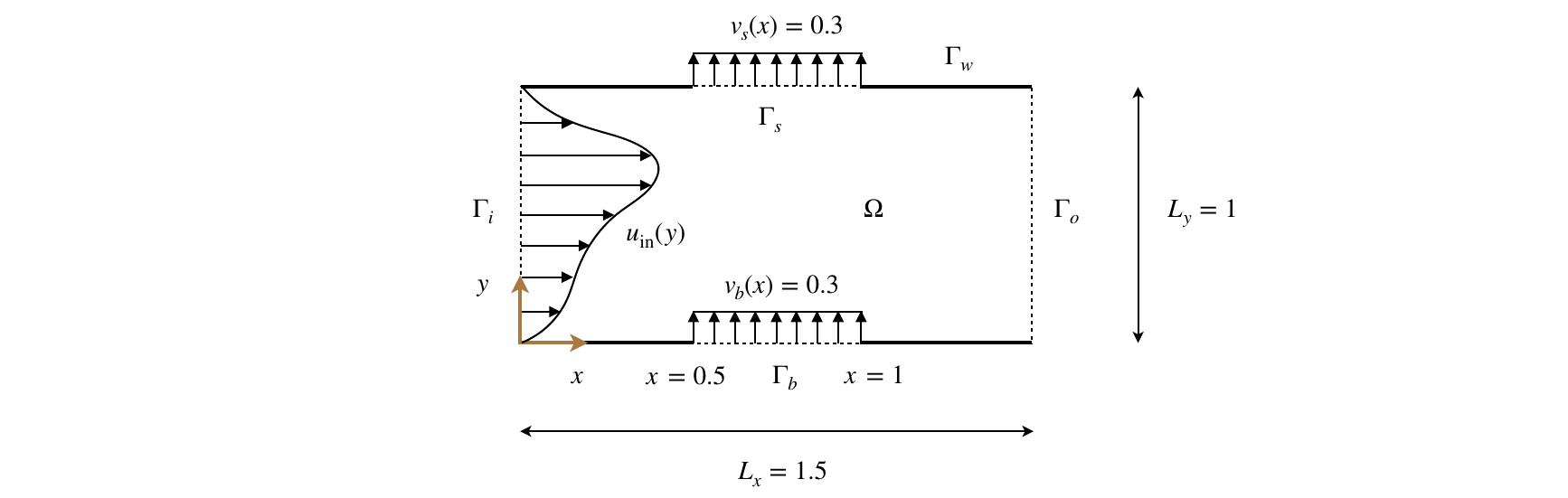}
\caption{\textcolor{black}{Scheme} of the setup for the Navier-Stokes equation. A prescribed horizontal velocity profile $\mathbf{u}(y) = (u_\mathrm{in}(y),0)$ is assigned to the inlet $\Gamma_i$. The blowing and suction boundaries $\Gamma_b$ and $\Gamma_s$ are assigned a uniform vertical velocity $\mathbf{u}(x) = (0.3,0)$. The flow leaves the domain at the outflow boundary $\Gamma_o$, and the remaining boundaries $\Gamma_w$ are no-slip walls.}
\label{fig:NSSchema}
\end{figure}

As a last example, we consider the steady 2D incompressible Navier-Stokes (NS) equations in the geometry depicted in Figure \ref{fig:NSSchema}.
In non-dimensional form, these equations are expressed as
\begin{subequations}
\begin{align}
(\mathbf{u} \cdot \nabla) \mathbf{u} &= - \nabla p + \frac{1}{Re} \nabla^2 \mathbf{u}, \label{eq:NSMom} \\
\nabla \cdot \mathbf{u} &= 0, \label{eq:NSCont}
\end{align} \label{eq:NS}%
\end{subequations}
where the velocity field $\mathbf{u}(\mathbf{x}) = (u(x,y),v(x,y))$ and pressure field $p(\mathbf{x}) = p(x,y)$ are defined in the rectangular 2D domain $\Omega = (L_x,L_y) = (1.5,1)$ shown in Figure \ref{fig:NSSchema}, and we choose the Reynolds number $Re = 100$. The boundary conditions for the velocity are
\begin{subequations}
\begin{align}
\mathbf{u} = (u_\mathrm{in}(y),0) &\quad \text{on } \Gamma_i, \\
\mathbf{u} = (v_b(x),0) &\quad \text{on } \Gamma_b, \\
\mathbf{u} = (v_s(x),0) &\quad \text{on } \Gamma_s, \\
(\mathbf{n} \cdot \nabla) \mathbf{u} = (0,0) &\quad \text{on } \Gamma_o, \\
\mathbf{u} = (0,0) &\quad \text{on } \Gamma_w,
\end{align}\label{eq:NSVelocityBCs}%
\end{subequations}
while the boundary conditions for the pressure are
\begin{subequations}
\begin{align}
(\mathbf{n} \cdot \nabla) p = 0 &\quad \text{on } \Gamma_i \cup \Gamma_b \cup \Gamma_s \cup \Gamma_w, \\
p = p_a &\quad \text{on } \Gamma_o,
\end{align} \label{eq:NSPressureBCs}%
\end{subequations}
where $\mathbf{n}$ denotes the unit surface normal, and $p_a$ is a reference pressure that we set to zero. These boundary conditions correspond to a prescribed horizontal velocity profile $u_\mathrm{in}(y)$ at an inlet $\Gamma_i$, a prescribed velocity profile $v_b(x)$ and $v_s(x)$ at two blowing and suction boundaries $\Gamma_b$ and $\Gamma_s$, an outflow boundary $\Gamma_o$ and no-slip walls $\Gamma_w$. We choose $v_b(x) = v_s(x) = 0.3$, and $u_\mathrm{in}(y)$ will be specified later. Although pressure boundary conditions are not usually stated explicitly, we will see below that they are important when evaluating the accuracy of the PINN solution.

\subsubsection{Forward problem}

We begin by solving a forward problem, defined by specifying a parabolic inlet velocity profile $u_\mathrm{in}(y) = u_\mathrm{parab}(y) = 4y(1-y)/L_y^2$. We compare the PINN solution to this problem with a reference numerical solution calculated using the finite-volume code OpenFOAM \cite{OpenFOAM}.

The PINN solution is obtained by representing $u(x,y)$, $v(x,y)$, $p(x,y)$ with a single network containing 5 hidden layers of 50 neurons each. The network takes $(x,y)$ as input and outputs $(u,v,p)$ at the corresponding location. It is trained by minimizing the loss \eqref{eq:Loss} using the Adam optimizer, starting from random initial weights. To evaluate the loss, we sample $N_r = 40000$ residual training points $(x_i,y_i) \in \Omega$ using LHS with 30000 points distributed in the entire domain and 10000 points distributed in 4 boxes of size $0.1 \times 0.02$ adjacent to the endpoints of $\Gamma_b$ and $\Gamma_s$. We select $N_b = 328$ boundary training points $(x_i,y_i)$, 82 of them equally spaced along the vertical boundaries $\Gamma_i$ and $\Gamma_o$, and the rest equally spaced along the horizontal boundaries $\Gamma_b$, $\Gamma_s$, and $\Gamma_w$. At the beginning of each epoch, the entire set of $N_r$ residual points is shuffled and divided into 10 minibatches of $N_r/10 = 4000$ points each. We set scalar weights $\lambda_{r}^{u\text{-mom}} = \lambda_{r}^{v\text{-mom}} = \lambda_{r}^\mathrm{cont} = 1$ and $w_b = 100$, with $\lambda_{r}^{u\text{-mom}}$, $\lambda_{r}^{v\text{-mom}}$ and $\lambda_{r}^\mathrm{cont}$ denoting the weights of the residuals corresponding to the $u$- and $v$-momentum equations \eqref{eq:NSMom} and continuity equation \eqref{eq:NSCont}. Finally, we choose an initial learning rate of $\alpha = 10^{-3}$ and decrease it by a factor 10 after 40k and 80k epochs, for a total of 100k epochs.

The finite-volume solution is carried out using the icoFOAM solver in OpenFOAM. Pressure and velocity are decoupled using the SIMPLE algorithm \cite{patankar1972calculation} technique. For the convection terms, second order Gaussian integration is used with the Sweby limiter \cite{sweby1984high} for numerical stability. For diffusion, Gaussian integration with central-differencing interpolation is used. The discretized algebraic equations are solved using the preconditioned biconjugate gradient (PBiCG) method. We employ a mesh of size 400 elements $(20\times 20)$, above which the cost function becomes independent of the mesh size. Anticipating \textcolor{black}{the DAL solution of the optimal control problem, for which we will need to compute the gradient of the cost objective with respect to the inlet velocity profile}, we perform local \textcolor{black}{mesh} refinements at the inlet for a better accuracy of the gradient calculation. Further, we found that using an upwind and first order method for solving the adjoint equations resulted in inaccurate gradients, and hence those methods were not adopted \cite{nabi2017adjoint,nabi2019nonlinear}.  

\begin{figure}
\centering
\includegraphics[width=\textwidth]{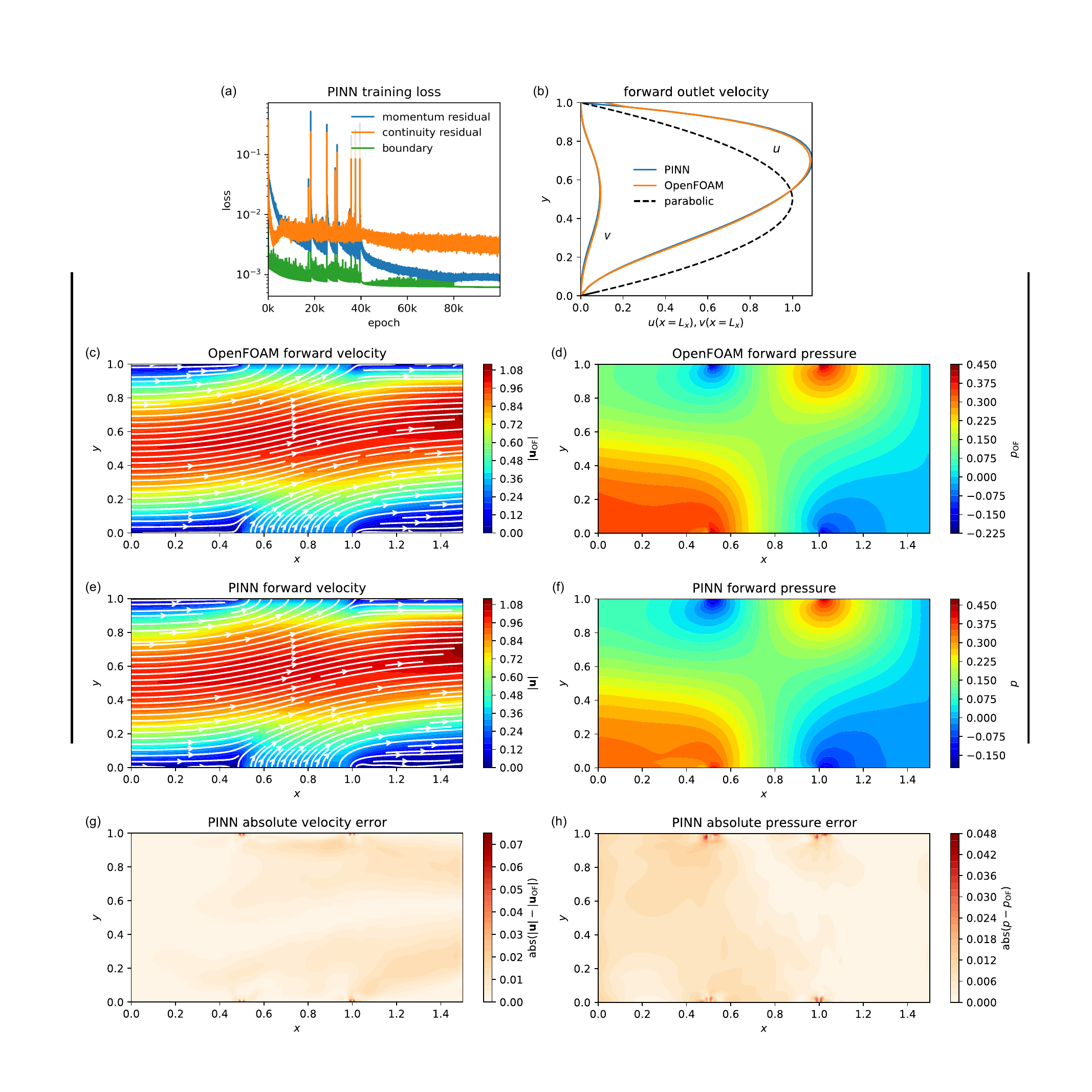}
\caption{Forward solution of the Navier-Stokes equation \eqref{eq:NS}. (a) Convergence of the loss during training of the PINN solution. (b) Velocity profile at the outlet $\Gamma_o$ of the trained PINN and the reference OpenFOAM solution. (c,e,g) Velocity magnitude and streamlines of the OpenFOAM and trained PINN solutions, and absolute error between them. (d,f,h) Pressure field of the OpenFOAM and trained PINN solutions.}
\label{fig:NSForward}
\end{figure}

The different loss components during the training of the PINN forward solution are shown in Figure \ref{fig:NSForward}(a).
The horizontal and vertical velocity profiles at the outlet boundary $\Gamma_o$ predicted by the trained PINN and reference OpenFOAM solutions are shown in Figure \ref{fig:NSForward}(b), displaying very good agreement with each other. The horizontal velocity profile is noticeably skewed upward, reaching its maximum value at $y = 0.7$. The velocity magnitude, streamlines and pressure fields corresponding to the OpenFOAM and PINN solutions are shown in Figures \ref{fig:NSForward}(c,d,e,f), and reveal that the skewed outlet velocity profile is caused by the blowing and suction boundaries, which deflect the parabolic inlet profile as the fluid moves through the channel. The low absolute error between the velocity magnitudes and pressure fields from the two solutions, shown in Figures \ref{fig:NSForward}(g,h), validates the choice of neural network parameters and training points.

We note that the literature on PINNs often mentions the nonnecessity of imposing pressure boundary conditions (BCs) in the PINN framework. While this is clearly advantageous when they are unknown or lack a clear physical motivation, pressure BCs do affect the resulting velocity profile. For our purpose of validating the PINN solution against a reference OpenFOAM solution, we thus found it important to include in the PINN loss function the same pressure BCs \eqref{eq:NSPressureBCs} employed in the OpenFOAM solver. In this way, we ensure that the PINN and OpenFOAM solvers are compared based on the exact same system of equations.

\subsubsection{Optimal control problem}

We observed in the forward problem that the parabolic inlet velocity profile $u_\mathrm{in}(y) = u_\mathrm{parab}(y)$ produces a skewed velocity profile at the outlet, Figure \ref{fig:NSForward}(b). Is it then possible to find an inlet velocity profile $u_\mathrm{in}^*(y)$ so that the outlet velocity profile is close to parabolic? This motivates the control problem
\begin{equation}
u_\mathrm{in}^* = \arg \min_{u_\mathrm{in}} \mathcal{J}(\mathbf{u}) \quad \text{subject to } \eqref{eq:NS}, \eqref{eq:NSVelocityBCs} \text{ and } \eqref{eq:NSPressureBCs},
\label{eq:OptimalControlNS}
\end{equation}
where the objective cost is
\begin{equation}
\mathcal{J}(\mathbf{u}) = \frac{1}{2} \int_0^{L_y} (|u(L_x,y) - u_\mathrm{parab}(y)|^2 + |v(L_x,y)|^2) dy, \quad u_\mathrm{parab}(y) = \frac{4}{L_y^2}y(1-y).
\label{eq:CostNS}
\end{equation}

To solve this problem in the PINNs framework, we define a second neural network for the control (inlet profile) $u_i(y)$, consisting of 3 hidden layers of 30 neurons each, in addition to the network for $u,v,p$ for which we keep the same architecture as before. We then train both networks simultaneously using the loss \eqref{eq:LossControl}, starting from a new initialization of the parameters. We use the same training points as before, and evaluate the integral in the cost objective \eqref{eq:CostKS} using the midpoint rule at $N_\mathcal{J} = 41$ equally-spaced training points on the outflow boundary $\Gamma_o$. We use the scalar weights $\lambda_{r}^{u\text{-mom}} = \lambda_{r}^{v\text{-mom}} = \lambda_{r}^\mathrm{cont} = 1$, $w_b = 100$. We choose an initial learning rate of $10^{-3}$ and decrease it by a factor 10 after 100k and 200k epochs of training, for a total of 300k epochs. \textcolor{black}{We repeat this procedure for 12 values of $w_\mathcal{J}$ between $10^{-3}$ and $10^5$.}

We implement the DAL procedure with the continuous adjoint formulation in the icoFoam solver of OpenFOAM, based on the adjoint NS equations and the gradient of the cost objective given in \ref{sec:AdjointNSEquation}. The adjoint equations are solved with the same numerical methods as the direct equations. The DAL optimal solution is obtained by iteratively solving the NS equations and their adjoint, updating the control $u_\mathrm{in}(y)$ at each iteration with the gradient descent formula \eqref{eq:GradientUpdateDAL}. We choose the parabolic velocity profile $u_\mathrm{parab}(y)$ as initial guess for the control $u_\mathrm{in}(y)$ and employ a learning rate $\beta = 0.001$.

\textcolor{black}{Beginning with Step 1 of the line search strategy presented in Section \ref{sec:Guidelines}, we visualize in Figure \ref{fig:NSControl}(a) the two components $\mathcal{L}_{\mathcal{F}/\mathcal{B}}$ and $\mathcal{L}_{\mathcal{J}}$ of the loss \eqref{eq:LossControlSimp} obtained at the end of training of the PINN optimal control solution, for each considered value of $w_\mathcal{J}$.
\begin{figure}
\centering
\includegraphics[width=\textwidth]{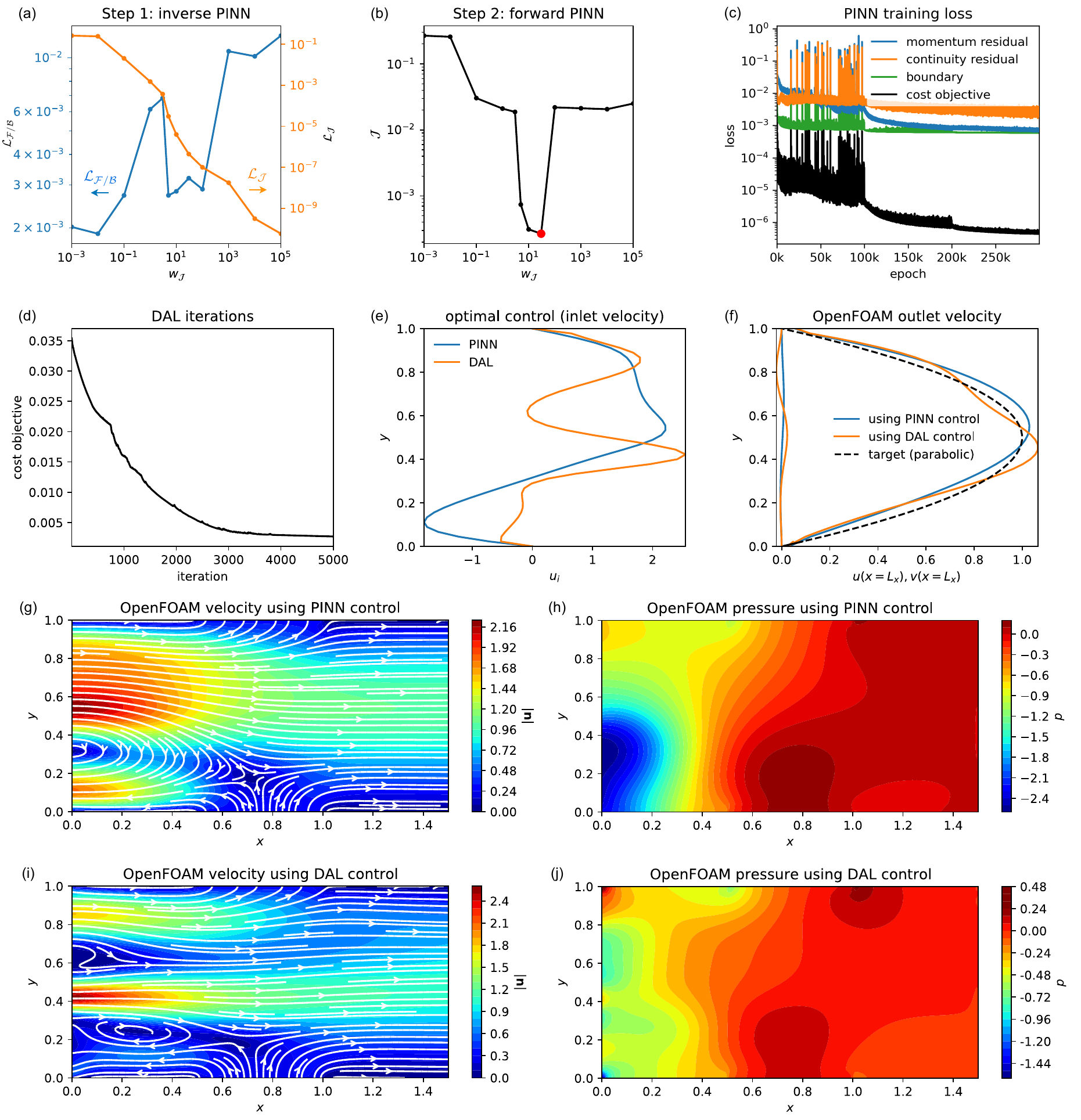}
\caption{Optimal solution of the Navier-Stokes control problem \eqref{eq:OptimalControlNS}. \textcolor{black}{(a) Components of the loss \eqref{eq:LossControlSimp} obtained at the end of training of the PINN optimal control solution versus weight $w_\mathcal{J}$ (step 1 of the line search strategy). (b) Cost objective estimate obtained by a separate PINN solution of the forward problem with fixed control from the PINN optimal solution versus $w_\mathcal{J}$ (step 2 of the line search strategy). The best optimal control, obtained with $w_\mathcal{J} = 30$, is shown by the red dot.} (c) Convergence of the loss during training of the PINN \textcolor{black}{optimal control solution (for $w_\mathcal{J} = 30$)}. (d) Convergence of the cost objective during DAL iterations. (e) Optimal inlet velocity profiles $u_\mathrm{in}^*$ obtained using the PINN \textcolor{black}{(for $w_\mathcal{J} = 30$)} and DAL frameworks. (f) Outlet velocity profiles of two forward OpenFOAM solutions calculated using the optimal inlet profiles $u_\mathrm{in}^*$ from the PINN \textcolor{black}{(for $w_\mathcal{J} = 30$)} and DAL frameworks, compared with the target parabolic profile. (f,g,h,i) Velocity magnitude, streamlines and pressure fields of the OpenFOAM solutions calculated using the optimal inlet profiles $u_\mathrm{in}^*$ from the PINN \textcolor{black}{(for $w_\mathcal{J} = 30$)} and DAL frameworks.}
\label{fig:NSControl}
\end{figure}
Moving to Step 2 of the line search strategy, we report in Figure \ref{fig:NSControl}(b) the cost objective value $\mathcal{J}$ obtained by a separate PINN solution of the corresponding forward problem, with fixed control from the optimal solution obtained in Step 1, for each value of $w_\mathcal{J}$. We observe that the control solution obtained with $w_\mathcal{J} = 30$, shown by the red dot, yields the lowest $\mathcal{J}$. It is therefore this optimal control solution that we analyze hereafter and compare with the DAL solution.}

The various loss components during training of the PINN optimal solution \textcolor{black}{(for $w_\mathcal{J} = 30$)} are displayed in Figure \ref{fig:NSControl}(c), and the convergence of the cost objective during the DAL iterations is shown in Figure \ref{fig:NSControl}(d). The PINN and DAL frameworks converge to rather different optimal inlet velocity profiles $u_\mathrm{in}^*(y)$, shown in Figure \ref{fig:NSControl}(e). These profiles nonetheless share a few features: two local maxima near the centerline and around $y = 0.8$, as well as a region of negative velocity for $y < 0.3$. We evaluate the quality of these optimal inlet profiles by computing the corresponding flow fields in two separate forward OpenFOAM calculations. The resulting outlet velocity profiles are displayed in Figures \ref{fig:NSControl}(f) and compared with the target parabolic profile $u_\mathrm{parab}(y)$. Both optimal inlet profiles lead to near-parabolic outlet profiles with comparable cost values: $\mathcal{J} = 0.00278$ and $\mathcal{J} = 0.00265$ for the PINN and DAL inlet profiles, respectively. Yet, the PINN inlet velocity profile is smoother and produces an outlet profile that has a more parabolic shape than its DAL counterpart. Figures \ref{fig:NSControl}(g,h,i,j) display the velocity magnitude, streamlines and pressure fields of the two OpenFOAM solutions calculated using the PINN and DAL optimal inlet profiles. In both cases, the bottom region of negative inlet velocity attracts some of the fluid entering through the blowing boundary, which reduces the effect of the latter on the outlet profile. The DAL inlet profile has a second region of negative velocity around $y=0.6$ and a much sharper peak near $y=0.4$, but these features do not yield a noticeably lower cost objective.

\section{Discussion}
\label{sec:Discussion}

A major goal of the present study is to compare the pros and cons of the PINN and DAL frameworks for solving PDE-constrained optimal control problems, so that the novel PINN approach can be placed in the context of the mature field of PDE-constrained optimization. To this end, in the previous section we have systematically compared the efficacy of the two frameworks by looking at several examples with varying degrees of complexity. In all cases, both techniques returned optimal control solutions yielding a comparable objective cost when evaluated in a separate high-fidelity numerical solver, which required an extra step in the PINN approach. The optimal control problem based on Laplace's equation had a global analytical optimal solution, which was found by both methods. The three remaining examples had more complicated optimization landscapes due to the nonlinearity of their governing PDEs, which exposed some differences in the characteristics of the optimal control solutions found by PINNs and DAL. For the Burgers and Navier-Stokes equations, the optimal control distributions found by DAL yielded a lower cost objective but were less smooth than the ones obtained from PINNs. For the Kuramoto-Sivashinsky equation, whose optimal control problem was the classical scenario of driving the state towards an unstable fixed point while penalizing large control magnitude, the two frameworks found remarkably close optimal control solutions yielding almost identical cost objective values. Such agreement is significant given the engineering relevance of this control problem and the chaotic nature of the Kuramoto-Sivashinsky equation.

Thus, these four examples suggest that the PINN approach can be similarly effective as DAL in solving optimal control problems. Clearly, an important advantage of the PINN framework is its ease of implementation, since it takes very little effort to adapt a PINN code for forward problems to the solution of optimal control problems. Furthermore, the PINN framework is very flexible in terms of the type of governing equations, boundary conditions, geometries, and cost objective functions that it allows, which opens the door for a very large class of problems to be solved with little human effort. By contrast, the DAL approach involves the cumbersome manual derivation of the adjoint equation and cost objective gradient, which needs to be repeated after any mere change of boundary conditions or cost objective function. The adjoint equation and DAL iterative procedure then need to be implemented in a numerical solver, which is no small task when complicated governing PDEs and/or complex geometries are involved. There exists a discrete adjoint formulation that does not require deriving adjoint equations, as opposed to the continuous adjoint formulation that we have employed in this paper. However, the discrete adjoint approach necessitates the use of differentiable numerical solvers, which are not yet well established. Thus, the PINN framework brings optimal control problems within reach of a much wider audience compared with adjoint-based approaches such as DAL.

We close the discussion with a few remarks regarding the computational cost of the PINN and DAL approaches for solving optimal control problems. The PINNs literature often alludes to the potentially superior computational efficiency of PINNs for solving inverse problems compared to other approaches, but actual numbers are rarely given. Indeed, it is very difficult to do an apples-to-apples comparison of the two frameworks due to the very different nature of their hardware requirements: PINNs are trained much faster on GPUs while the forward and adjoint equations in DAL are solved on CPUs. In addition, the computational cost of PINNs depends on many factors such as the machine learning library that is used, the optimization algorithm, and the neural network architecture. Likewise, the computational cost of DAL depends on other factors such as the discretization method and programming language employed to solve the forward and adjoint equations as well as the gradient descent algorithm. For instance, BFGS-type updates can be used to speed up the convergence of DAL.
With that in mind, our goal is not to do a rigorous comparison of the PINN and DAL computational costs, but rather to provide the reader with an idea of what to expect, which is still valuable given the lack of such data in the literature. To this effect, we report in Table \ref{tab:ComputationalCost} the computational times required to obtain all the PINN and DAL optimal control solutions reported in the previous section.
\begin{table*}
\begin{center}
\renewcommand{\arraystretch}{1.1}
\begin{tabular}{|r|cc|}
\hline
 & PINN & DAL
 \\ \hline
Laplace (Section \ref{sec:Laplace}) & 9 min (Tesla V100) & 25 min (Xeon E5-2683) \\
Burgers (Section \ref{sec:Burgers}) & 19 min (Tesla V100) & 1 min (Core i7-4980HQ) \\
Kuramoto-Sivashinsky (Section \ref{sec:KS}) & 2 hours 4 min (Tesla V100) & 2 hours 55 min (Core i7-4980HQ) \\
Navier-Stokes (Section \ref{sec:NS}) & 8 hours 20 min (Tesla V100) & 28 hours (Xeon E5-2683) \\\hline
\end{tabular}
\caption{Computational time for obtaining the PINN and DAL optimal solutions. The PINN solutions are trained on one GPU in TensorFlow. The DAL solutions are all computed using a single CPU core, using the C++ finite-volume solver OpenFOAM for the Laplace and Navier-Stokes equations, and a spectral Python code for the Burgers and Kuramoto-Sivashinsky equations. Note that this table is not meant as a rigorous comparison of the computational efficiency of the two frameworks due to the number of differing factors involved.}
\label{tab:ComputationalCost}
\end{center}
\end{table*}
For the simpler problems based on the Laplace and Burgers equations, the DAL solution was obtained in much shorter time than the PINN solution\footnote{Note that the DAL iterations for the Laplace solution were pushed to an extremely low value of the cost objective; a quarter of the iterations would have produced an optimal control of a similar quality to the PINN solution.}. For the more complex problems based on the Kuramoto-Sivashinsky and Navier-Stokes equations, the situation reverses and the PINN solution is obtained in shorter time. This suggests that the computational efficiency of the PINN framework might close the gap with that of DAL as the problem  complexity increases. There are, however, a few caveats to this picture to keep in mind.
Importantly, all our DAL calculations are performed on a single CPU core, and could therefore benefit from a substantial speedup if performed in a multi-threaded environment. Further, although both frameworks necessitate some amount of parameter tuning, the process is more involved for PINNs due to the number of parameters involved, even when following the guidelines presented in Section \ref{sec:Guidelines}.
Finally, algorithms other than DAL to solve the adjoint-based optimality equations exist, e.g.~one-shot methods, that may result in overall computational cost of optimization comparable to that of a single forward simulation \cite{borzi2011computational,bosse2014one}. 

\section{Conclusions}
\label{sec:Conclusions}

In this paper, we have proposed a methodology and a set of guidelines for solving optimal control problems with physics-informed neural networks (PINNs). We then compared rigorously optimal control solutions obtained from PINNs with corresponding results calculated with direct-adjoint-looping (DAL), a particular implementation of adjoint-based optimization which is the standard approach for solving PDE-constrained optimal control problems. The comparison was carried out over four examples with increasing complexity levels: the Laplace, Burgers, Kuramoto-Sivashinsky, and Navier-Stokes equations. In all cases, both techniques found optimal control solutions that yielded comparable cost objective values once plugged back into a high-fidelity numerical solver, demonstrating the capability of PINNs to solve optimal control problems. Beyond the quality of the optimal control solutions, we also assessed the pros and cons of the two approaches. A major strength of the PINNs framework is its flexibility and ease of implementation, which hold the potential to make optimal control problems more accessible to a much wider audience than adjoint-based methods. On the other hand, DAL is very time-consuming to carry out but can potentially return the optimal control solution at a smaller computational cost if the solver is parallelized. This being said, the performance of the PINN framework may also be improved by using recent advances in the field such as adaptive weighting of the various loss components \cite{wang2020understanding,maddu2022inverse} or adaptive refinement of the residual points \cite{lu2021deepxde}. Finally, PINNs will automatically benefit from the future advances of fast-evolving deep-learning tools, so their relevance for optimal control problems will only increase from now on.

\appendix

\section{Adjoint Burgers equation}
\label{sec:AdjointBurgersEquation}

Applying the methodology outlined in Section \ref{sec:AdjointBasedOptimalControl}, one can derive the adjoint of the Burgers equation \eqref{eq:Burgers} as well as the gradient of the cost objective \eqref{eq:CostBurgers} with respect to the control. The adjoint Burgers equation, obtained from \eqref{eq:AdjointPDE} using integration by part, is
\begin{equation}
- \frac{\partial \lambda}{\partial t} - u \frac{\partial \lambda}{\partial x} = \nu \frac{\partial^2 \lambda}{\partial x^2}, 
\end{equation}
where $\lambda(x,t)$ is the adjoint field and $u(x,t)$ is the forward field that solves the Burgers equation \eqref{eq:Burgers} given the control (initial condition) $u_0(x)$. The adjoint equation is supplemented with periodic boundary conditions and the terminal condition $\lambda(x,T) = -u + u_d$. Finally, the total gradient of the cost objective with respect to the control $u_0(x)$, obtained from \eqref{eq:CostGradient}, is
\begin{equation}
\frac{\mathrm{d}\mathcal{J}(u)}{\mathrm{d}u_0} = -\lambda(x,0).
\end{equation}

\section{Adjoint Kuramoto-Sivashinsky equation}
\label{sec:AdjointKSEquation}

Applying the methodology outlined in Section \ref{sec:AdjointBasedOptimalControl}, one can derive the adjoint of the KS equation \eqref{eq:KS} as well as the gradient of the cost objective \eqref{eq:CostKS} with respect to the control. The adjoint KS equation, obtained from \eqref{eq:AdjointPDE} using integration by part, is
\begin{equation}
- \frac{\partial \lambda}{\partial t} - u \frac{\partial \lambda}{\partial x} + \frac{\partial^2 \lambda}{\partial x^2} + \frac{\partial^4 \lambda}{\partial x^4} =  -u(x,t),
\end{equation}
where $\lambda(x,t)$ is the adjoint field and $u(x,t)$ is the forward field that solves the KS equation \eqref{eq:KS} given the control (forcing) $f(x,t)$. The adjoint equation is supplemented with periodic boundary conditions and the terminal condition $\lambda(x,T) = 0$. Finally, the total gradient of the cost objective with respect to the control $f(x,t)$, obtained from \eqref{eq:CostGradient}, is
\begin{equation}
\frac{\mathrm{d}\mathcal{J}(u,f)}{\mathrm{d}f} = 2 \sigma f(x,t) - \lambda(x,t).
\end{equation}

\section{Adjoint Navier-Stokes equations}
\label{sec:AdjointNSEquation}

In this appendix, we use Einstein notation so that the velocity field defined in Section \ref{sec:NS} will be denoted $\mathbf{u}(\mathbf{x}) = (u_1(x_1,x_2),u_2(x_1,x_2))$. The augmented objective functional, i.e.~Lagrangian, corresponding to the control problem \eqref{eq:OptimalControlNS} is
\begin{equation}\label{app1}
\mathcal{L}=\mathcal{J}+ \bigg\langle \lambda_i,\frac{\partial {u}_i{u}_j}{\partial x_j} + \frac{\partial {p}_i}{\partial x_i}-\frac{1}{Re}\frac{\partial^2 u_i}{\partial x_j^2} \bigg\rangle + \bigg\langle \Pi,-\frac{\partial u_j}{\partial x_j} \bigg\rangle,
\end{equation} 
where $\boldsymbol{\lambda} = (\lambda_1(x_1,x_2),\lambda_2(x_1,x_2))$ and $\Pi$ are adjoint velocity and pressure fields, respectively, and the inner product $\langle \cdot, \cdot \rangle$ is defined as
\begin{equation}
\langle a, b \rangle = \int_\Omega a(\mathbf{x}) b(\mathbf{x}) dV.
\end{equation}
The variation of the Lagrangian is
\begin{equation}\label{app2}
\delta \mathcal{L}=\delta \mathcal{J} + \bigg\langle \lambda_i,\delta{u}_j\frac{\partial {u}_i}{\partial x_j}+{u}_j\frac{\partial \delta u_i}{\partial x_j}+\frac{\partial \delta p_i}{\partial x_i}- \frac{1}{Re}\frac{\partial^2 \delta u_j}{\partial x_j^2} \bigg\rangle + \bigg\langle \Pi,-\frac{\partial \delta u_j}{\partial x_j} \bigg\rangle.
\end{equation}
For optimality $\delta \mathcal{L}=0$ should be satisfied. Using vector calculus and integration by parts for each term, appropriate Euler-Lagrange equations can be derived. 
For instance, 
\begin{subequations}
\begin{align}
\bigg\langle \lambda_i,\delta{u}_j\frac{\partial {u}_i}{\partial x_j} \bigg\rangle &= \bigg\langle \delta u_i,\lambda_j\frac{\partial {u}_j}{\partial x_i} \bigg\rangle, \\
\bigg\langle \lambda_i,{u}_j\frac{\partial \delta u_i}{\partial x_j} \bigg\rangle &= -\bigg\langle \delta u_i,u_j\frac{\partial \lambda_i}{\partial x_j} \bigg\rangle + \int_{\partial \Omega} \lambda_i\delta u_i u_jn_j dS, \\
\bigg\langle \lambda_i,\frac{1}{Re}\frac{\partial^2 \delta u_i}{\partial x_j^2} \bigg\rangle &= \bigg\langle \delta u_i,\frac{1}{Re}\frac{\partial^2 \lambda_i}{\partial x_j^2} \bigg\rangle + \int_{\partial \Omega} \frac{1}{Re} \left( n_j\frac{\partial \lambda_i}{\partial x_j}\delta u_i-n_j\frac{\partial \delta u_i}{\partial x_j} \lambda_i \right) dS,
\end{align} \label{eq_app6}%
\end{subequations}
and so on for the other terms. Here, $\mathbf{n} = (n_1,n_2)$ is the normal unit vector of the surface. From the volumetric integrals, the adjoint equations are recovered as
\begin{subequations}
\begin{align}
\lambda_j\frac{\partial u_j}{\partial x_i}-u_j\frac{\partial \lambda_i}{\partial x_j}-\frac{1}{Re}\frac{\partial^2 \lambda_i}{\partial x_j^2}+\frac{\partial \Pi}{\partial x_i}&=0, \\
\frac{\partial \lambda_j}{\partial x_j}&=0.
\end{align} \label{eq_app7}%
\end{subequations}
Setting surface integrals to zero and decomposing these integrals into normal and tangential components, we obtain the corresponding boundary conditions for the adjoint velocity and pressure as
\begin{subequations}
\begin{align}
\lambda_1 = \lambda_2 = 0, (\mathbf{n} \cdot \nabla) \Pi = 0 &\quad \text{on } \Gamma_i \cup \Gamma_b \cup \Gamma_s \cup \Gamma_w, \\
u_1 \lambda_2 + \frac{1}{Re} \frac{\partial \lambda_2}{\partial x_1} = -u_2 &\quad \text{on } \Gamma_o, \\
u_1 \lambda_1 + \frac{1}{Re} \frac{\partial \lambda_1}{\partial x_1} = u_1-u_\mathrm{parab} + \Pi &\quad \text{on } \Gamma_o.
\end{align}\label{eq:eq_app8}%
\end{subequations}
Finally, the total gradient of the cost objective with respect to the control, or design equation, is given by 
\begin{equation}
\frac{\mathrm{d}\mathcal{J}(\mathbf{u})}{\mathrm{d}u_\mathrm{in}} =\Pi(0,x_2) - \frac{1}{Re} \frac{\partial \lambda_1}{\partial x_1}(0,x_2),
\label{eq_app9}
\end{equation}
where all values are evaluated at the inlet. For more details, the reader is invited to refer to \citep{nabi2019nonlinear}.

\bibliography{bibliography}

\end{document}